\documentclass[12pt]{article}
\usepackage[utf8]{inputenc}
\usepackage{mathrsfs}
\usepackage{fullpage}
\usepackage{amssymb}
\usepackage{amsmath}
\usepackage{epsfig}
\usepackage{algorithm}
\usepackage{algorithmic}
\usepackage{latexsym}
\usepackage{amsmath}
\usepackage{amsfonts}
\usepackage{graphicx}
\usepackage{authblk}
\usepackage[normalem]{ulem}
\usepackage{array}
\usepackage{colortbl}
\usepackage{bm}
\usepackage{color}\usepackage[dvipsnames]{xcolor}
\usepackage{tikz}\usetikzlibrary{patterns}
\usepackage{chessboard}
\usepackage{skak}
\usepackage{array,multirow}
\usepackage[
  bookmarks=false, 
  colorlinks,
  citecolor=brown!70!black,
  linkcolor=brown!80!black,
  urlcolor=blue!70!black,
]{hyperref}
\usepackage{relsize,exscale}

\newtheorem{cor}{Corollary}
\newtheorem{thm}{Theorem}
\newtheorem{rem}{Remark}

\newtheorem{prop}{Proposition}

\newcommand{\comm}[1]{}
 	\definecolor{lightlightgray}{rgb}{0.93, 0.93, 0.93}
 		\definecolor{llightgray}{rgb}{0.87, 0.87, 0.87}

\newcolumntype{C}[1]{>{\centering\arraybackslash }b{#1}}

\def\A{0.5cm}

 \def\C{2.5cm}
 
 \def\E{4.5cm}

 \def\G{6.5cm}
 
 \def\I{8.5cm}

 \def\K{10.5cm}
 
\def\M{12.5cm}
 
 \def\O{14.5cm}\def\Q{16.5cm}\def\S{18.5cm}
\def\U{20.5cm}\def\W{22.5cm}\def\Y{24.5cm}\def\ZZ{26.5cm}\def\Za{28.5cm}\def\Zb{30.5cm}\def\Zc{32.5cm}\def\Zd{34.5cm}\def\Ze{36.5cm}

\newcommand\scalemath[2]{\scalebox{#1}{\mbox{\ensuremath{\displaystyle #2}}}}

\newcommand\oeis[1]{\href{https://oeis.org/#1}{#1}}

\DeclareMathAlphabet{\mymathbb}{U}{BOONDOX-ds}{m}{n}

\definecolor{bleuclair}{RGB}{186,183,229}
\definecolor{rougeclair}{RGB}{255,230,231}
\definecolor{bleufonce}{RGB}{59,50,114}
\definecolor{rougefonce}{RGB}{127,0,4}

\let\ge\geqslant

\let\geq\geqslant

\title{Two kinds of partial Motzkin paths with air pockets}

\author[1]{Jean-Luc Baril}
\author[2]{Paul Barry}
\affil[1]{\rm LIB, Universit\'e de Bourgogne Franche-Comt\'e \protect\\
  B.P. 47 870, 21078 Dijon Cedex France\protect\\
   {\tt E-mail: barjl@u-bourgogne.fr
   }
}

\affil[2]{\rm School of Science, South East Technological University (SETU)\protect\\ Ireland\protect\\
   {\tt E-mail: pbarry@wit.ie
   }
}
\date{\today}

\date{\today}
\begin{document}

\maketitle

\begin{abstract}
 Motzkin paths with air pockets (MAP) are defined as a generalization of Dyck paths with air pockets by adding some horizontal steps with certain conditions. In this paper, we introduce two  generalizations. The first one consists of lattice paths in $\Bbb{N}^2$ starting at the origin made of steps $U=(1,1)$, $D_k=(1,-k)$, $k\geq 1$ and $H=(1,0)$, where two down steps cannot be consecutive, while the second one are lattice paths in $\Bbb{N}^2$ starting at the origin, made of steps $U$, $D_k$ and $H$, where each step $D_k$ and $H$ is necessarily followed by an up step, except for the last step of the path. We provide enumerative results for these paths according to the length, the type of the last step, and the height of its end-point. A similar study is made for these paths read from right to left. As a byproduct, we obtain new classes of paths counted by the  Motzkin numbers. Finally, we express our results using Riordan arrays.
\end{abstract}

\section{Introduction}
In a recent  paper \cite{bakimava}, the authors introduce, study and enumerate special classes of lattice paths, called {\it Dyck paths with air pockets} (DAP for short). Such paths are  non empty lattice paths in the first quadrant of $\Bbb{Z}^2$ starting at the origin, and consisting of up-steps $U=(1,1)$ and down-steps $D_k=(1,-k)$, $k\geq 1$, where two down steps cannot be consecutive. These paths can be viewed as ordinary Dyck paths where each maximal run of down-steps is condensed into one large down step. As mentioned in \cite{bakimava}, they also correspond to a stack evolution with (partial) reset operations that cannot be consecutive (see for instance \cite{krin}). The authors enumerate these paths with respect to the length, the type (up or down) of the last step and the height of the end-point. Whenever the last point is on the $x$-axis, they prove that the DAP of length $n$ are in one-to-one correspondence with the peakless Motzkin paths of length $n-1$.  They also investigate the popularity of many patterns in these paths and they give asymptotic approximations. In a second work \cite{bkmv}, the authors make a study for a generalization of these paths by allowing them to go below the $x$-axis. They call these paths Grand Dyck paths with air pockets (GDAP), and they also yield enumerative results for these paths according to the length and several restrictions on the height.

In this paper, we  introduce two generalizations of partial Dyck paths of air pockets by allowing some possible horizontal steps $H=(1,0)$ with some conditions. These two kinds of paths can be viewed as special partial  Motzkin paths (lattice paths in $\Bbb{N}^2$ starting at the origin and made of $U$, $D$, and $H$), where each maximal run of down-steps is condensed into one large down step.

Firstly, we consider lattice paths in $\Bbb{N}^2$ starting at the origin, consisting of steps $U$, $D_k$ and $H$, where {\it two down steps cannot be consecutive}. Secondly, we consider lattices paths in $\Bbb{N}^2$ starting at the origin, consisting of steps $U$, $D_k$ and $H$, {\it where any step $U$ and $D_k$ (except the last step of the path) is immediately followed by an up step $U$}. These two classes of paths will be denoted $\mathcal{M}_1$ and $\mathcal{M}_2$, respectively. The paths in $\mathcal{M}_1$ (resp. $\mathcal{M}_2$) will be called {\it partial Motzkin paths with air pockets} of first kind (resp. second kind), and they are called {\it Motzkin paths with air pockets} whenever they end on the $x$-axis. For short, we denote by PMAP  all  paths in $\mathcal{M}_i$, $i\in\{1,2\}$, and by PAP all paths ending on the $x$-axis.
On the other hand, let $\mathcal{M}'_1$ (resp. $\mathcal{M}'_2$) be the set of lattice paths starting at the origin obtained by reading the paths in $\mathcal{M}_1$ (resp. $\mathcal{M}_2$) from right to left, i.e. up steps are changed into down step and {\it vice versa}, and horizontal steps are unchanged (see below for a more formal definition of these paths).

Throughout the paper, and for each class of paths $\mathcal{M}_i$ and $\mathcal{M}'_i$, $i\in\{1,2\}$, described above,  we will use the following notations. For $k\geq 0$, we consider the generating function $f_k=f_k(z)$ (resp. $g_k=g_k(z)$, resp. $h_k=h_k(z)$), where the coefficient of $z^n$ in the series expansion is the number of  partial Motzkin paths with air pockets of length $n$  ending at height $k$ with an up-step, (resp. with a down-step,   resp. with a horizontal step $H$). 

We introduce the bivariate generating functions
$$F(u,z)=\sum\limits_{k\geq 0} u^kf_k(z), \quad G(u,z)=\sum\limits_{k\geq 0} u^kg_k(z), \mbox{ and } H(u,z)=\sum\limits_{k\geq 0} u^kh_k(z).$$
For short, we also use the notation $F(u), G(u)$ and $H(u)$ for these functions.
\medskip

The outline of this paper is the following. In Section 2, we present enumerative results for partial Motzkin paths with air pockets of the first kind, and for these paths when we read them from right to left. We provide bivariate generating functions that count these paths with respect to the length, the type of the last step (up, down or horizontal step) and the height of the end-point. In Section 3, we make a similar study for PMAP of second kind, and we present new classes of lattice paths counted by the well known Motzkin numbers. All these results are obtained algebraically by using the famous kernel method for solving several systems of functional equations. Finally, we express our results using Riordan arrays and we deduce closed forms for  PMAP of length $n$ ending at height $k$.

\section{PMAP of the first kind }
In this section, we focus on  PMAP of the first kind, i.e. lattices paths in $\Bbb{N}^2$ starting at the origin, made of steps $U$, $D_k$ and $H$, such that two down steps cannot be consecutive. The first subsection considers the paths  in $\mathcal{M}_1$, while the second subsection handles the paths in $\mathcal{M}'_1$ (see Introduction for the definition of these two sets). We yield enumerative results for these paths according to the length, the type of the last step, and the height of its end-point. 

\subsection{PMAP in $\mathcal{M}_1$ - From left to right }
\begin{figure}[h]
 \begin{center}
        \begin{tikzpicture}[scale=0.15]
            \draw (\A,\A)-- (38,\A);
             \draw[dashed,line width=0.1mm] (\A,\E)-- (\Ze,\E);
              \draw[dashed,line width=0.1mm] (\A,\C)-- (\Ze,\C);
               \draw[dashed,line width=0.1mm] (\A,\G)-- (\Ze,\G);
               \draw[dashed,line width=0.1mm] (\A,\I)-- (\Ze,\I);
              \draw[dashed,line width=0.1mm] (\A,\K)-- (\Ze,\K);
               \draw[dashed,line width=0.1mm] (\A,\M)-- (\Ze,\M);
            \draw (\A,\A) -- (\A,\O);
             \draw[dashed,line width=0.1mm] (\C,\A) -- (\C,\M);\draw[dashed,line width=0.1mm] (\E,\A) -- (\E,\M);\draw[dashed,line width=0.1mm] (\G,\A) -- (\G,\M);
             \draw[dashed,line width=0.1mm] (\I,\A) -- (\I,\M);\draw[dashed,line width=0.1mm] (\K,\A) -- (\K,\M);\draw[dashed,line width=0.1mm] (\M,\A) -- (\M,\M);
             \draw[dashed,line width=0.1mm] (\O,\A) -- (\O,\M);\draw[dashed,line width=0.1mm] (\Q,\A) -- (\Q,\M);\draw[dashed,line width=0.1mm] (\S,\A) -- (\S,\M);
             \draw[dashed,line width=0.1mm] (\U,\A) -- (\U,\M);\draw[dashed,line width=0.1mm] (\W,\A) -- (\W,\M);\draw[dashed,line width=0.1mm] (\Y,\A) -- (\Y,\M);
             \draw[dashed,line width=0.1mm] (\ZZ,\A) -- (\ZZ,\M);
             \draw[dashed,line width=0.1mm] (\Za,\A) -- (\Za,\M);
             \draw[dashed,line width=0.1mm] (\Zb,\A) -- (\Zb,\M);
             \draw[dashed,line width=0.1mm] (\Zc,\A) -- (\Zc,\M);
             \draw[dashed,line width=0.1mm] (\Zd,\A) -- (\Zd,\M);
             \draw[dashed,line width=0.1mm] (\Ze,\A) -- (\Ze,\M);
            \draw[solid,line width=0.4mm] (\A,\A)--(\C,\C)  -- (\E,\C) -- (\G,\C) --(\I,\A)-- (\K,\C) -- (\M,\E) -- (\O,\E) -- (\Q,\G)  -- (\S,\A)--(\W,\E)--(\Y,\E)--(\Za,\I) -- (\Zb,\C) -- (\Zd, \G)--(\Ze,\A);

         \end{tikzpicture}\qquad 
 \begin{tikzpicture}[scale=0.15]
            \draw (\A,\A)-- (38,\A);
             \draw[dashed,line width=0.1mm] (\A,\E)-- (\Ze,\E);
              \draw[dashed,line width=0.1mm] (\A,\C)-- (\Ze,\C);
               \draw[dashed,line width=0.1mm] (\A,\G)-- (\Ze,\G);
               \draw[dashed,line width=0.1mm] (\A,\I)-- (\Ze,\I);
              \draw[dashed,line width=0.1mm] (\A,\K)-- (\Ze,\K);
               \draw[dashed,line width=0.1mm] (\A,\M)-- (\Ze,\M);
            \draw (\A,\A) -- (\A,\O);
             \draw[dashed,line width=0.1mm] (\C,\A) -- (\C,\M);\draw[dashed,line width=0.1mm] (\E,\A) -- (\E,\M);\draw[dashed,line width=0.1mm] (\G,\A) -- (\G,\M);
             \draw[dashed,line width=0.1mm] (\I,\A) -- (\I,\M);\draw[dashed,line width=0.1mm] (\K,\A) -- (\K,\M);\draw[dashed,line width=0.1mm] (\M,\A) -- (\M,\M);
             \draw[dashed,line width=0.1mm] (\O,\A) -- (\O,\M);\draw[dashed,line width=0.1mm] (\Q,\A) -- (\Q,\M);\draw[dashed,line width=0.1mm] (\S,\A) -- (\S,\M);
             \draw[dashed,line width=0.1mm] (\U,\A) -- (\U,\M);\draw[dashed,line width=0.1mm] (\W,\A) -- (\W,\M);\draw[dashed,line width=0.1mm] (\Y,\A) -- (\Y,\M);
             \draw[dashed,line width=0.1mm] (\ZZ,\A) -- (\ZZ,\M);
             \draw[dashed,line width=0.1mm] (\Za,\A) -- (\Za,\M);
             \draw[dashed,line width=0.1mm] (\Zb,\A) -- (\Zb,\M);
             \draw[dashed,line width=0.1mm] (\Zc,\A) -- (\Zc,\M);
             \draw[dashed,line width=0.1mm] (\Zd,\A) -- (\Zd,\M);
             \draw[dashed,line width=0.1mm] (\Ze,\A) -- (\Ze,\M);
            \draw[solid,line width=0.4mm] (\A,\A)--(\C,\C)  -- (\E,\C) -- (\G,\C) --(\I,\A)-- (\K,\C) -- (\M,\E) -- (\O,\E) -- (\Q,\G)  -- (\S,\A)--(\W,\E)--(\Y,\E)--(\Za,\I) -- (\Zb,\C) -- (\Zc, \E)--(\Zd,\E)--(\Ze,\G);

         \end{tikzpicture}
               \end{center}
         \caption{ The left drawing shows a Motzkin path with air pockets of length $18$. The right drawing shows a partial Motzkin path with air pockets of  length $18$  ending at height $3$.}
         \label{fig1}
\end{figure}
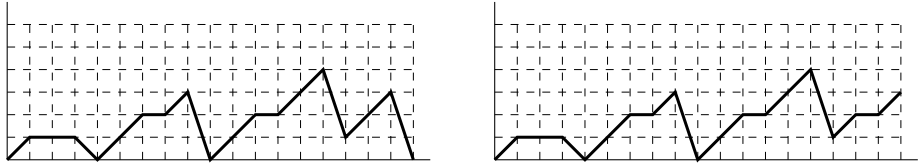

 In this part, we consider PMAP in $\mathcal{M}_1$. Figure \ref{fig1} shows two examples of such paths.

 Let $P$ be a length $n$ PMAP in $\mathcal{M}_1$ ending at height $k\geq 0$. 
 If the last step of $P$ is $U$, then $k\geq 1$ and we have $P=QU$ where $Q$ is a length $(n-1)$ MPAP ending at height $k-1$. So, we obtain the first relation  $f_k=zf_{k-1}+zg_{k-1}+zh_{k-1}$ for $k\geq 1$, anchored with $f_0=1$ by considering the empty path.  If the last step of $P$ is a down step $D_k$, $k\geq 0$, then we have $P=QD_k$ where $Q$ is a length $(n-1)$ PMAP ending at height $\ell\geq k+1$ with no down step at its end. So, we obtain the second relation  $g_k=z\sum\limits_{\ell\geq k+1} f_\ell+z\sum\limits_{\ell\geq k+1} h_\ell$. If the last step of $P$ is a horizontal step $H$,  then we have $P=QH$ where $Q$ is a length $(n-1)$ PMAP ending at height $k$. 
 
 Therefore, we have to solve the following system of equations.
 
\begin{equation}\left\{\begin{array}{l}
f_0=1,\mbox{ and } f_k=zf_{k-1}+zg_{k-1}+zh_{k-1}, \quad k\geq 1,\\
g_k=z\sum\limits_{\ell\geq k+1} f_\ell+z\sum\limits_{\ell\geq k+1} h_\ell,\quad k\geq 0,\\
h_k=zf_k+zg_k+zh_k,\quad k\geq 0.\\
\end{array}\right.\label{equ1}
\end{equation}

Summing the recursions in (\ref{equ1}), we have:
\begin{align*}
F(u)&=1+z\sum\limits_{k\geq 1}u^kf_{k-1}+z\sum\limits_{k\geq 1}u^kg_{k-1} +z\sum\limits_{k\geq 1}u^kh_{k-1}\\
&=1+zuF(u)+zuG(u)+zuH(u),\\
G(u)&=z\sum\limits_{k\geq 0}u^k \Bigl(\sum\limits_{\ell\geq k+1} f_\ell     \Bigr)+z\sum\limits_{k\geq 0}u^k \Bigl(\sum\limits_{\ell\geq k+1} h_\ell     \Bigr)\\
&=z\sum\limits_{k\geq 1}f_k(1+u+\ldots + u^{k-1})+ z\sum\limits_{k\geq 1}h_k(1+u+\ldots + u^{k-1})\\
&=z\sum\limits_{k\geq 1}\frac{u^k-1}{u-1}f_k+z\sum\limits_{k\geq 1}\frac{u^k-1}{u-1}h_k\\
&=\frac{z}{u-1}(F(u)-F(1)+H(u)-H(1)),\\
H(u)&=zF(u)+zG(u)+zH(u).
\end{align*}

Notice that we have $F(1)-H(1)=1$ by considering the difference of the first and third equations. Now, setting $f_1:=F(1)$ and solving these functional equations, we obtain
$$F(u)=\frac{2 f_1 u \,z^{2}-u \,z^{2}+u z +z^{2}-u -z +1}{u^{2} z +u \,z^{2}+z^{2}-u -z +1},$$
$$G(u)=-\frac{z \left(2 f_1 u z +2 zf_1 -u z -2 f_1 -z +2\right)}{u^{2} z +u \,z^{2}+z^{2}-u -z +1},\quad 
H(u)=\frac{z \left(2 z f_1 -u -2 z +1\right)}{u^{2} z +u \,z^{2}+z^{2}-u -z +1}.$$

In order to compute $f_1$, we use the kernel method (see~\cite{ban, pro}) on $F(u)$. We can write the denominator (which is a polynomial in $u$ of degree 2), as $z(u-r)(u-s)$ with 
$$r=\frac{-z^{2}+1+\sqrt{z^{4}-4 z^{3}+2 z^{2}-4 z +1}}{2 z}, \mbox{ and }
s=\frac{1-z^{2}-\sqrt{z^{4}-4 z^{3}+2 z^{2}-4 z +1}}{2 z}.$$
Plugging $u=s$ (which have a Taylor expansion at $z=0$) in $F(u)z(u-r)(u-s)$, we obtain the equation 
$$2 f_1 s\,z^{2}-s \,z^{2}+s z +z^{2}-s -z +1 =0.$$
 Using $zrs=z^2-z+1$, we deduce
 $$f_1=F(1)=H(1)+1=\frac{r(s-1)}{2z}.$$
Finally,  after simplifying by the factor $(u-s)$ in the numerators and denominators, we obtain 
 $$F(u)=\frac{r}{r-u},\quad G(u)=\frac{r(s-1)-z}{r-u},\quad\mbox{ and } \quad H(u)=\frac{1}{r-u},$$
which implies that
$$f_k=[u^k] F(u)=\frac{1}{r^k}, \quad g_k=[u^k] G(u)=\frac{s-1}{r^k}-\frac{z}{r^{k+1}},\quad \mbox{ and }\quad
h_k=[u^k] H(u)=\frac{1}{r^{k+1}}.
$$

\begin{thm} The bivariate generating function for the total number of PMAP in $\mathcal{M}_1$ with respect to the  length  and the height of the end-point is given by
$$\mathit{Total}(z,u)=\frac{1+rs-z}{r-u},$$ and we have
$$[u^k] \mathit{Total}(z,u)=\frac{1+rs-z}{r^{k+1}}.$$

Finally, setting $t(n,k)=[z^n][u^k]\mathit{Total}(z,u)$, we have for $n\geq 2$, $k\geq 1$,  $$t(n,k)=t(n,k-1)+t(n-1,k)-t(n-1,k-2)-t(n-2,k)-t(n-2,k-1),
$$
and setting $t_n:=t(n,0)$, then we have $$t_n=t_{n-1}+t_{n-2}+\sum\limits_{k=0}^{n-3}t_{k}t_{n-k-3}+\sum\limits_{k=2}^{n-1}\left(t_k-t_{k-1}\right)t_{n-k-1}.$$
\end{thm}
\noindent {\it Proof.} 
The first two equalities are immediately deduced from the previous results. The third equality is obtained by checking that 
$$\mathit{Total}(z,u)=(u+z-zu^2-z^2-z^2u)\mathit{Total}(z,u)-u+\frac{(z^2-z+1)(1+rs-z)}{r}.$$
Now, let us prove the last equality. Any length $n$ MAP is of the form ($i$) $HP$, or ($ii$) $UDP$, or ($iii$) $U P FD Q$ where $P,Q$ are some MAP so that the length of $P$ lies into $[0,n-3]$, or ($iv$) $P^\sharp Q$  where $P^\sharp=UP'D_{k}$, $k\geq 1$, and $P'D_{k-1}$ is a MAP  of length lying into $[2,n-1]$. Taking into account all these cases, we obtain the result.  
\hfill $\Box$

Let $\mathcal{T}$ be the infinite matrix $\mathcal{T}:=[t(n,k)]_{n\geq 0, k\geq 0}$. The first few rows of the matrix $\mathcal{T}$ are
\[\mathcal{T}=
\left(
\begin{array}{ccccccccc}
 1 & 0 & 0 & 0 & 0 & 0 & 0 & 0 & 0 \\
 1 & 1 & 0 & 0 & 0 & 0 & 0 & 0 & 0 \\
 2 & 2 & 1 & 0 & 0 & 0 & 0 & 0 & 0 \\
 5 & 5 & 3 & 1 & 0 & 0 & 0 & 0 & 0 \\
 13 & 14 & 9 & 4 & 1 & 0 & 0 & 0 & 0 \\
 36 & 40 & 28 & 14 & 5 & 1 & 0 & 0 & 0 \\
 105 & 118 & 87 & 48 & 20 & 6 & 1 & 0 & 0 \\
 317 & 359 & 273 & 161 & 75 & 27 & 7 & 1 & 0 \\
 982 & 1118 & 869 & 536 & 270 & 110 & 35 & 8 & 1 \\
\end{array}
\right).\]

\begin{cor} The g.f. that counts the PMAP  in $\mathcal{M}_1$ with respect to the length is given by 
$$\mathit{Total}(z,1)=\frac{1+rs-z}{r-1}.$$
\end{cor}

The first few terms of the series expansion of $\mathit{Total}(z,1)$ are  $1 + 2z + 5z^2 + 14z^3 + 41z^4 + 124z^5 + 385z^6 + 1220z^7 + 3929z^8 + 12822z^9+ O(z^{10})$, which correspond to the sequence \oeis{A159771} in \cite{oeis} counting the $n$-leaf binary trees that do not contain $(()((()())((()())())))$ as a subtree (see \cite{row}).

\begin{cor} The g.f. that counts the MAP in  $\mathcal{M}_1$ with respect to the length is given by 
$$\mathit{Total}(z,0)=\frac{1+rs-z}{r}.$$
\end{cor}

The first few terms of the series expansion of $\mathit{Total}(z,0)$ are
$1 + z + 2z^2 + 5z^3 + 13z^4 + 36z^5 + 105z^6 + 317z^7 + 982z^8 + 3105z^9 + O(x^{10})$
which correspond to the sequence \oeis{A114465} in \cite{oeis} counting  Dyck paths of semilength $n$ having no ascents of length $2$ that start at an odd level. We leave open the question of finding a constructive bijection between these sets.

\subsection{PMAP in $\mathcal{M}'_1$ - From right to left}
\begin{figure}[h]
 \begin{center}
        \begin{tikzpicture}[scale=0.15]
            \draw (\A,\A)-- (38,\A);
             \draw[dashed,line width=0.1mm] (\A,\E)-- (\Ze,\E);
              \draw[dashed,line width=0.1mm] (\A,\C)-- (\Ze,\C);
               \draw[dashed,line width=0.1mm] (\A,\G)-- (\Ze,\G);
               \draw[dashed,line width=0.1mm] (\A,\I)-- (\Ze,\I);
              \draw[dashed,line width=0.1mm] (\A,\K)-- (\Ze,\K);
               \draw[dashed,line width=0.1mm] (\A,\M)-- (\Ze,\M);
            \draw (\A,\A) -- (\A,\O);
             \draw[dashed,line width=0.1mm] (\C,\A) -- (\C,\M);\draw[dashed,line width=0.1mm] (\E,\A) -- (\E,\M);\draw[dashed,line width=0.1mm] (\G,\A) -- (\G,\M);
             \draw[dashed,line width=0.1mm] (\I,\A) -- (\I,\M);\draw[dashed,line width=0.1mm] (\K,\A) -- (\K,\M);\draw[dashed,line width=0.1mm] (\M,\A) -- (\M,\M);
             \draw[dashed,line width=0.1mm] (\O,\A) -- (\O,\M);\draw[dashed,line width=0.1mm] (\Q,\A) -- (\Q,\M);\draw[dashed,line width=0.1mm] (\S,\A) -- (\S,\M);
             \draw[dashed,line width=0.1mm] (\U,\A) -- (\U,\M);\draw[dashed,line width=0.1mm] (\W,\A) -- (\W,\M);\draw[dashed,line width=0.1mm] (\Y,\A) -- (\Y,\M);
             \draw[dashed,line width=0.1mm] (\ZZ,\A) -- (\ZZ,\M);
             \draw[dashed,line width=0.1mm] (\Za,\A) -- (\Za,\M);
             \draw[dashed,line width=0.1mm] (\Zb,\A) -- (\Zb,\M);
             \draw[dashed,line width=0.1mm] (\Zc,\A) -- (\Zc,\M);
             \draw[dashed,line width=0.1mm] (\Zd,\A) -- (\Zd,\M);
             \draw[dashed,line width=0.1mm] (\Ze,\A) -- (\Ze,\M);
            \draw[solid,line width=0.4mm] (\A,\A)--(\C,\G)  -- (\E,\E) -- (\G,\C) --(\I,\I)-- (\K,\G) -- (\M,\E) -- (\O,\E) -- (\Q,\I)  -- (\S,\I)--(\U,\G)--(\W,\E)--(\Y,\E)--(\Za,\A) -- (\Zb,\G) --(\Ze,\A);

         \end{tikzpicture}\qquad 
 \begin{tikzpicture}[scale=0.15]
            \draw (\A,\A)-- (38,\A);
             \draw[dashed,line width=0.1mm] (\A,\E)-- (\Ze,\E);
              \draw[dashed,line width=0.1mm] (\A,\C)-- (\Ze,\C);
               \draw[dashed,line width=0.1mm] (\A,\G)-- (\Ze,\G);
               \draw[dashed,line width=0.1mm] (\A,\I)-- (\Ze,\I);
              \draw[dashed,line width=0.1mm] (\A,\K)-- (\Ze,\K);
               \draw[dashed,line width=0.1mm] (\A,\M)-- (\Ze,\M);
            \draw (\A,\A) -- (\A,\O);
             \draw[dashed,line width=0.1mm] (\C,\A) -- (\C,\M);\draw[dashed,line width=0.1mm] (\E,\A) -- (\E,\M);\draw[dashed,line width=0.1mm] (\G,\A) -- (\G,\M);
             \draw[dashed,line width=0.1mm] (\I,\A) -- (\I,\M);\draw[dashed,line width=0.1mm] (\K,\A) -- (\K,\M);\draw[dashed,line width=0.1mm] (\M,\A) -- (\M,\M);
             \draw[dashed,line width=0.1mm] (\O,\A) -- (\O,\M);\draw[dashed,line width=0.1mm] (\Q,\A) -- (\Q,\M);\draw[dashed,line width=0.1mm] (\S,\A) -- (\S,\M);
             \draw[dashed,line width=0.1mm] (\U,\A) -- (\U,\M);\draw[dashed,line width=0.1mm] (\W,\A) -- (\W,\M);\draw[dashed,line width=0.1mm] (\Y,\A) -- (\Y,\M);
             \draw[dashed,line width=0.1mm] (\ZZ,\A) -- (\ZZ,\M);
             \draw[dashed,line width=0.1mm] (\Za,\A) -- (\Za,\M);
             \draw[dashed,line width=0.1mm] (\Zb,\A) -- (\Zb,\M);
             \draw[dashed,line width=0.1mm] (\Zc,\A) -- (\Zc,\M);
             \draw[dashed,line width=0.1mm] (\Zd,\A) -- (\Zd,\M);
             \draw[dashed,line width=0.1mm] (\Ze,\A) -- (\Ze,\M);
            \draw[solid,line width=0.4mm](\A,\A)--(\C,\G)  -- (\E,\E) -- (\G,\C) --(\I,\I)-- (\K,\G) -- (\M,\E) -- (\O,\E) -- (\Q,\I)  -- (\S,\I)--(\U,\G)--(\W,\E)--(\Y,\E)--(\Za,\A) -- (\Zb,\G) --(\Zc,\E)--(\Zd,\G)--(\Ze,\E);;

         \end{tikzpicture}
               \end{center}
         \caption{ The left drawing shows a Motzkin path with air pockets of length $18$ read from right to left. The right drawing shows a partial Motzkin path with air pockets of  length $18$  ending at height $2$ and reaf from right to left.}
         \label{fig2}
\end{figure}
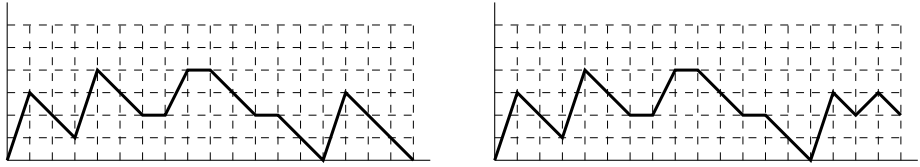

Here, we consider the paths of the previous section, but we read them from right to left. This means that down steps become up steps  and {\it vice versa}, and horizontal steps are unchanged, which implies that two up steps cannot be consecutive now. See Figure \ref{fig2} for two examples of such paths.

With the same arguments and the same notations used in the previous part,  we can easily obtain the following equations:
\begin{equation}\left\{\begin{array}{l}
f_0=1,\mbox{ and } f_k=z+z(g_0+g_1+\ldots +g_{k-1})+z(h_0+h_1+\ldots +h_{k-1}), \quad k\geq 1,\\
g_k=zf_{k+1}+zg_{k+1}+zh_{k+1},\quad k\geq 0,\\
h_k=zf_k+zg_k+zh_k,\quad k\geq 0.\\
\end{array}\right.\label{equ2}
\end{equation}

Summing the recursions in (\ref{equ2}), we have:
\begin{align*}
F(u)&=1+\frac{zu}{1-u}+z\sum\limits_{k\geq 1}(g_0+\ldots +g_{k-1})u^k+z\sum\limits_{k\geq 1}(h_0+\ldots +h_{k-1})u^k\\
&=1+\frac{zu}{1-u}+z\sum\limits_{k\geq 0}g_k\frac{u^{k+1}}{1-u}+z\sum\limits_{k\geq 0}h_k\frac{u^{k+1}}{1-u}\\
&=1+\frac{zu}{1-u}(1+G(u)+H(u)),\\
G(u)&=z\sum\limits_{k\geq 0}f_{k+1}u^k+z\sum\limits_{k\geq 0}g_{k+1}u^k+z\sum\limits_{k\geq 0}h_{k+1}u^k\\
&=\frac{z}{u}(F(u)-F(0)+G(u)-G(0)+H(u)-H(0)),\\
H(u)&=zF(u)+zG(u)+zH(u).
\end{align*}

Notice that we have $H(0)=\frac{z(1+G(0))}{1-z}$ by considering the third equation. Now, setting $g_0:=G(0)$ and solving these functional equations, we obtain
$$F(u)={\frac {-{u}^{2}{z}^{3}+g_0\,u{z}^{2}+3\,{u}^{2}{z}^{2}-u{z}^{3}
-3\,{u}^{2}z+2\,u{z}^{2}+{u}^{2}+uz-{z}^{2}-u+z}{ \left( 1-z \right) 
 \left( {u}^{2}{z}^{2}-{u}^{2}z+u{z}^{2}+{u}^{2}-u+z \right) }}
,$$
$$G(u)={\frac {z \left( g_0\,u{z}^{2}-g_0\,uz+g_0\,u+g_0
\,z+uz-g_0 \right) }{ \left( -1+z
 \right)  \left( {u}^{2}{z}^{2}-{u}^{2}z+u{z}^{2}+{u}^{2}-u+z \right)}}
,$$
$$H(u)={\frac { \left( g_0\,u{z}^{2}+{u}^{2}{z}^{2}-g_0\,uz-2\,{u}
^{2}z+u{z}^{2}+g_0\,z+{u}^{2}-u+z \right) z}{ \left( 1-z
 \right)  \left( {u}^{2}{z}^{2}-{u}^{2}z+u{z}^{2}+{u}^{2}-u+z \right) 
}}
.$$

In order to compute $g_0$, we use the kernel method (see~\cite{ban, pro}) on $G(u)$. We can write the denominator (which is a polynomial in $u$ of degree 2), as $(z-1)(z^2-z+1)(u-r)(u-s)$ with 
$$r={\frac {1-{z}^{2}+\sqrt {{z}^{4}-4\,{z}^{3}+2\,{z}^{2}-4\,z+1}}{
2({z}^{2}-z+1)}}
, \mbox{ and }
s={\frac {1-{z}^{2}-\sqrt {{z}^{4}-4\,{z}^{3}+2\,{z}^{2}-4\,z+1}}{
2({z}^{2}-z+1)}}
.$$
Plugging $u=s$ (which have a Taylor expansion at $z=0$) in $G(u)(z^2-z+1)(u-r)(u-s)$, we obtain the equation 
$$g_0(sz^2-sz+s+z-1)+sz=0.$$
Using $sr(z^2-z+1)=z$, we deduce 
 $$g_0=G(0)=\frac{1-r}{r}-sz.$$

Finally,  after simplifying by the factor $(u-s)$ in the numerators and denominators, we obtain

$$F(u)=1+\frac{sru}{r-u},\quad G(u)=\frac{s(1-r+rz)}{r-u},\quad\mbox{ and } \quad H(u)=\frac{sr-sru(1-z)}{r-u},$$

which implies that
\begin{align}f_k&=[u^k] F(u)=[k=0]+(1-[k=0])\cdot\frac{s}{r^{k-1}},\\
g_k&=[u^k] G(u)=\frac{s(1-r+rz)}{r^{k+1}},\\
h_k&=[u^k] H(u)=\frac{s}{r^{k}}-(1-[k=0])\cdot\frac{(1-z)s}{r^{k-1}}.
\end{align}

\begin{thm} The bivariate generating function for the total number of PMAP (read from right to left)  with respect to the length and the height of the end-point is given by
$$\mathit{Total}(z,u)=1+\frac{s(1+rz+ruz)}{r-u},$$ and we have
$$[u^k] \mathit{Total}(z,u)=[k=0]+\frac{s(rz+1)}{r^{k+1}}+(1-[k=0])\cdot\frac{sz}{r^{k-1}}.$$

Finally, setting $t(n,k)=[z^n][u^k]\mathit{Total}(z,u)$, we have for $n\geq 2$, $k\geq 1$, 
$$t(n,k)=t(n-2,  k-1) + t(n-2,  k) - t(n-1,  k-1) + t(n-1, k+1) + t(n,  k-1),
$$
and setting $t_n:=t(n,0)$, we have $$t_n=t_{n-1}+t_{n-2}+\sum\limits_{k=0}^{n-3}t_{k}t_{n-k-3}+\sum\limits_{k=2}^{n-1}\left(t_k-t_{k-1}\right)t_{n-k-1}.$$
\end{thm}
\noindent {\it Proof.} The first two equalities are directly deduced from the previous results. Since  the expression
$$(u-u^2z^2-uz^2+zu^2-z-u^2)\mathit{Total}(z,u)+u(1-z)$$ is a polynomial of degree one in $u$, we deduce the third relation. The last equality is already given in Theorem 1.
\hfill $\Box$

Let $\mathcal{T}$ be the infinite matrix $\mathcal{T}:=[t(n,k)]_{n\geq 0, k\geq 0}$. The first few rows of the matrix $\mathcal{T}$ are
\[\mathcal{T}=
\left(
\begin{array}{ccccccc}
 1 & 0 & 0 & 0 & 0 & 0 & 0  \\
 1 & 1 & 1 & 1 & 1 & 1 & 1  \\
 2 & 3 & 3 & 3 & 3 & 3 & 3  \\
 5 & 8 & 10 & 12 & 14 & 16 & 18  \\
 13 & 23 & 33 & 43 & 53 & 63 & 73  \\
 36 & 69 & 107 & 149 & 195 & 245 & 299  \\
 105 & 212 & 348 & 512 & 704 & 924 & 1172  \\
 317 & 665 & 1141 & 1753 & 2509 & 3417 & 4485  \\
 982 & 2123 & 3771 & 5999 & 8879 & 12483 & 16883 \\
\end{array}
\right).\]
Since there is an infinite number of PMAP of length $n$, we do not provide an ordinary generating function (with
respect to the length) for these paths. So, we get around this  by counting PMAP ending on a point $(x,n-x)$ for a given $n\geq 0$.

\begin{cor} The g.f. that counts the partial PMAP ending on the line $y=n-x$ is given by 
$$\mathit{Total}(z,z)=1+\frac{s(1+rz+rz^2)}{r-z}.$$
\end{cor}

The first few terms of the series expansion of $\mathit{Total}(z,z)$ are  $1 + z + 3z^2 + 9z^3 + 25z^4 + 73z^5 + 223z^6 + 697z^7 + 2217z^8 + 7161z^9+ O(z^{10})$, which correspond to the sequence \oeis{A101499} in \cite{oeis}, which is a Chebyshev transform of the Catalan number that counts peakless Motzkin paths of length $n$ where horizontal steps at level at least one come in $2$ colors.

Notice that we obviously retrieve (see Corollary 2) that the g.f. counting the MAP with respect to the length is given by 
$$\mathit{Total}(z,0)=1+\frac{s(1+rz)}{r}.$$


\section{PMAP of the second kind }

In this section, we focus on  PMAP of the second kind. The first subsection considers paths in $\mathcal{M}_2$, while  the second handles paths in $\mathcal{M}'_2$. We yield enumerative results for these paths according to the length, the type of the last step, and the height of the end-point.

\subsection{PMAP in $\mathcal{M}_2$ - From left to right}

\begin{figure}[h]
 \begin{center}
        \begin{tikzpicture}[scale=0.15]
            \draw (\A,\A)-- (38,\A);
             \draw[dashed,line width=0.1mm] (\A,\E)-- (\Ze,\E);
              \draw[dashed,line width=0.1mm] (\A,\C)-- (\Ze,\C);
               \draw[dashed,line width=0.1mm] (\A,\G)-- (\Ze,\G);
               \draw[dashed,line width=0.1mm] (\A,\I)-- (\Ze,\I);
              \draw[dashed,line width=0.1mm] (\A,\K)-- (\Ze,\K);
               \draw[dashed,line width=0.1mm] (\A,\M)-- (\Ze,\M);
            \draw (\A,\A) -- (\A,\O);
             \draw[dashed,line width=0.1mm] (\C,\A) -- (\C,\M);\draw[dashed,line width=0.1mm] (\E,\A) -- (\E,\M);\draw[dashed,line width=0.1mm] (\G,\A) -- (\G,\M);
             \draw[dashed,line width=0.1mm] (\I,\A) -- (\I,\M);\draw[dashed,line width=0.1mm] (\K,\A) -- (\K,\M);\draw[dashed,line width=0.1mm] (\M,\A) -- (\M,\M);
             \draw[dashed,line width=0.1mm] (\O,\A) -- (\O,\M);\draw[dashed,line width=0.1mm] (\Q,\A) -- (\Q,\M);\draw[dashed,line width=0.1mm] (\S,\A) -- (\S,\M);
             \draw[dashed,line width=0.1mm] (\U,\A) -- (\U,\M);\draw[dashed,line width=0.1mm] (\W,\A) -- (\W,\M);\draw[dashed,line width=0.1mm] (\Y,\A) -- (\Y,\M);
             \draw[dashed,line width=0.1mm] (\ZZ,\A) -- (\ZZ,\M);
             \draw[dashed,line width=0.1mm] (\Za,\A) -- (\Za,\M);
             \draw[dashed,line width=0.1mm] (\Zb,\A) -- (\Zb,\M);
             \draw[dashed,line width=0.1mm] (\Zc,\A) -- (\Zc,\M);
             \draw[dashed,line width=0.1mm] (\Zd,\A) -- (\Zd,\M);
             \draw[dashed,line width=0.1mm] (\Ze,\A) -- (\Ze,\M);
            \draw[solid,line width=0.4mm] (\A,\A)--(\C,\C)  -- (\E,\E) -- (\G,\G) --(\I,\A)-- (\K,\C) -- (\M,\E) -- (\O,\E) -- (\Q,\G)  -- (\S,\A)--(\W,\E)--(\Y,\E)--(\Za,\I) -- (\Zb,\C) -- (\Zd, \G)--(\Ze,\A);

         \end{tikzpicture}\qquad 
 \begin{tikzpicture}[scale=0.15]
            \draw (\A,\A)-- (38,\A);
             \draw[dashed,line width=0.1mm] (\A,\E)-- (\Ze,\E);
              \draw[dashed,line width=0.1mm] (\A,\C)-- (\Ze,\C);
               \draw[dashed,line width=0.1mm] (\A,\G)-- (\Ze,\G);
               \draw[dashed,line width=0.1mm] (\A,\I)-- (\Ze,\I);
              \draw[dashed,line width=0.1mm] (\A,\K)-- (\Ze,\K);
               \draw[dashed,line width=0.1mm] (\A,\M)-- (\Ze,\M);
            \draw (\A,\A) -- (\A,\O);
             \draw[dashed,line width=0.1mm] (\C,\A) -- (\C,\M);\draw[dashed,line width=0.1mm] (\E,\A) -- (\E,\M);\draw[dashed,line width=0.1mm] (\G,\A) -- (\G,\M);
             \draw[dashed,line width=0.1mm] (\I,\A) -- (\I,\M);\draw[dashed,line width=0.1mm] (\K,\A) -- (\K,\M);\draw[dashed,line width=0.1mm] (\M,\A) -- (\M,\M);
             \draw[dashed,line width=0.1mm] (\O,\A) -- (\O,\M);\draw[dashed,line width=0.1mm] (\Q,\A) -- (\Q,\M);\draw[dashed,line width=0.1mm] (\S,\A) -- (\S,\M);
             \draw[dashed,line width=0.1mm] (\U,\A) -- (\U,\M);\draw[dashed,line width=0.1mm] (\W,\A) -- (\W,\M);\draw[dashed,line width=0.1mm] (\Y,\A) -- (\Y,\M);
             \draw[dashed,line width=0.1mm] (\ZZ,\A) -- (\ZZ,\M);
             \draw[dashed,line width=0.1mm] (\Za,\A) -- (\Za,\M);
             \draw[dashed,line width=0.1mm] (\Zb,\A) -- (\Zb,\M);
             \draw[dashed,line width=0.1mm] (\Zc,\A) -- (\Zc,\M);
             \draw[dashed,line width=0.1mm] (\Zd,\A) -- (\Zd,\M);
             \draw[dashed,line width=0.1mm] (\Ze,\A) -- (\Ze,\M);
            \draw[solid,line width=0.4mm] (\A,\A)--(\C,\C)  -- (\E,\E) -- (\G,\G) --(\I,\A)-- (\K,\C) -- (\M,\E) -- (\O,\E) -- (\Q,\G)  -- (\S,\A)--(\W,\E)--(\Y,\E)--(\Za,\I) -- (\Zb,\C) -- (\Zc, \E)--(\Zd,\E)--(\Ze,\G);

         \end{tikzpicture}
               \end{center}
         \caption{ The left drawing shows a MAP of length $18$ in $\mathcal{M}_2$. The right drawing shows a PMAP of  length $18$ ending at height $3$ in $\mathcal{M}_2$.}
         \label{fig3}
\end{figure}
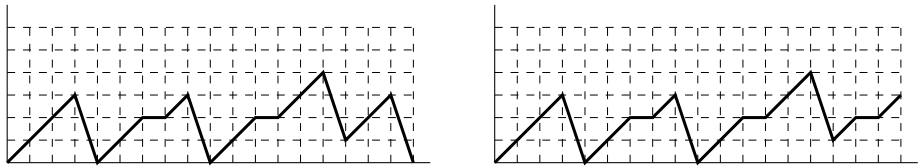

In this part, we consider PMAP in $\mathcal{M}_2$, i.e. lattice paths in $\Bbb{N}^2$ starting at the origin, consisting of steps $U$, $D_k$ and $H$, and where any down step or horizontal step  (except for the last step of the path) is immediately followed by an up step. Figure \ref{fig3} shows two examples of such paths.

 So, we easily obtain the following equations:
\begin{equation}\left\{\begin{array}{l}
f_0=1,\mbox{ and } f_k=zf_{k-1}+zg_{k-1}+zh_{k-1}, \quad k\geq 1,\\
g_k=z\sum\limits_{\ell\geq k+1} f_\ell,\quad k\geq 0,\\
h_k=zf_k,\quad k\geq 0.\\
\end{array}\right.\label{equ3}
\end{equation}

Summing the recursions in (\ref{equ3}), we have:
\begin{align*}
F(u)&=1+z\sum\limits_{k\geq 1}u^kf_{k-1}+z\sum\limits_{k\geq 1}u^kg_{k-1} +z\sum\limits_{k\geq 1}u^kh_{k-1}\\
&=1+zuF(u)+zuG(u)+zuH(u),\\
G(u)&=z\sum\limits_{k\geq 0}u^k \Bigl(\sum\limits_{\ell\geq k+1} f_\ell     \Bigr)=z\sum\limits_{k\geq 1}\frac{u^k-1}{u-1}f_k\\
&=\frac{z}{u-1}(F(u)-F(1)),\\
H(u)&=zF(u).
\end{align*}

Notice that $H(1)=zF(1)$ using the third equation. Now, setting $f_1:=F(1)$  and solving these functional equations, we deduce
$$F(u)=\frac{f_1 u \,z^{2}-u +1}{u^{2} z^{2}+u^{2} z -u z -u +1},$$
$$G(u)=-\frac{z \left(f_1 u \,z^{2}+f_1 u z -f_1 +1\right)}{u^{2} z^{2}+u^{2} z -u z -u +1},
H(u)=\frac{z \left(f_1 u \,z^{2}-u +1\right)}{u^{2} z^{2}+u^{2} z -u z -u +1}
.$$

In order to compute $f_1$, we use the kernel method (see~\cite{ban, pro}) on $F(u)$. We can write the denominator (which is a polynomial in $u$ of degree 2), as $z^2(u-r)(u-s)$ with 
$$r=\frac{z +1+\sqrt{-3 z^{2}-2 z +1}}{2 z \left(z +1\right)},\mbox{ and }
s=-\frac{-z -1+\sqrt{-3 z^{2}-2 z +1}}{2 z \left(z +1\right)}.$$
Plugging $u=s$ (which have a Taylor expansion at $z=0$) in $F(u)z^2(u-r)(u-s)$, we obtain the equation 
$$f_1 s\,z^{2}-s +1 =0,$$
and thus 
 $$f_1=\frac{s-1}{sz^2}.$$

Finally using $z(1+z)rs=1$ and simplifying by the factor $(u-s)$ in the numerators and denominators, we obtain

$$F(u)=\frac{r}{r-u},\quad
G(u)=\frac{s-1}{sz(r-u)}
,\quad \mbox{ and } \quad
H(u)=\frac{zr}{r-u}
,$$

which implies that
\begin{align}f_k&=[u^k] F(u)=\frac{1}{r^k},\\
g_k&=[u^k] G(u)=(1+z)\cdot\frac{s-1}{r^{k}},\\
h_k&=[u^k] H(u)=\frac{z}{r^k}.
\end{align}

\begin{thm} The bivariate generating function for the total number of PMAP with respect to the length and the height of the end-point is given by
$$\mathit{Total}(z,u)=\frac{1}{z(r-u)},$$ and we have
$$[u^k] \mathit{Total}(z,u)=\frac{1}{zr^{k+1}}.$$
Finally, setting $t(n,k)=[z^n][u^k]\mathit{Total}(z,u)$, we have for $n\geq 2$, $k\geq 1$,  $$t(n,k)=t(n,k-1)+t(n-1,k-1)-t(n-1,k-2)-t(n-2,k-2),
$$
and setting $t_n:=t(n,0)$, we have $$t_n=t_{n-1}+\sum\limits_{k=1}^{n-2}t_kt_{n-1-k}.$$
\end{thm}
\noindent {\it Proof.} The first two equalities are immediately deduced from the previous results. The third equality is obtained by checking that 
$$\mathit{Total}(z,u)=(u+zu-zu^2-z^2u^2)\mathit{Total}(z,u)-(1+z)u+\frac{1}{zr}.$$
For the last equality, it suffices to remark that the o.g.f of the first column, that is $1/(zr)$, generates a shift of the well known Motzkin sequence \oeis{A001006}. 
\hfill $\Box$

Let $\mathcal{T}$ be the infinite matrix $\mathcal{T}:=[t(n,k)]_{n\geq 0, k\geq 0}$. The first few rows of the matrix $\mathcal{T}$ are
\[\mathcal{T}=
\left(
\begin{array}{ccccccccc}
 1 & 0 & 0 & 0 & 0 & 0 & 0 & 0 & 0 \\
 1 & 1 & 0 & 0 & 0 & 0 & 0 & 0 & 0 \\
 1 & 2 & 1 & 0 & 0 & 0 & 0 & 0 & 0 \\
 2 & 3 & 3 & 1 & 0 & 0 & 0 & 0 & 0 \\
 4 & 6 & 6 & 4 & 1 & 0 & 0 & 0 & 0 \\
 9 & 13 & 13 & 10 & 5 & 1 & 0 & 0 & 0 \\
 21 & 30 & 30 & 24 & 15 & 6 & 1 & 0 & 0 \\
 51 & 72 & 72 & 59 & 40 & 21 & 7 & 1 & 0 \\
 127 & 178 & 178 & 148 & 105 & 62 & 28 & 8 & 1 \\
\end{array}
\right).\]

\begin{cor} The g.f. that counts the PMAP with respect to the length is given by 
$$\mathit{Total}(z,1)=\frac{1}{z(r-1)}.$$
\end{cor}

\begin{figure}[h]
 \begin{center}
        \begin{tikzpicture}[scale=0.15]
            \draw (\A,\A)-- (\I,\A);
           
             \draw[dashed,line width=0.1mm] (\A,\A)-- (\I,\A);
              \draw[dashed,line width=0.1mm] (\A,\C)-- (\I,\C);
               \draw[dashed,line width=0.1mm] (\A,\E)-- (\I,\E);
               \draw[dashed,line width=0.1mm] (\A,\G)-- (\I,\G);
               \draw[dashed,line width=0.1mm] (\A,\I)-- (\I,\I);
            \draw (\A,\A) -- (\A,\I);
             \draw[dashed,line width=0.1mm] (\C,\A) -- (\C,\I);\draw[dashed,line width=0.1mm] (\E,\A) -- (\E,\I);\draw[dashed,line width=0.1mm] (\G,\A) -- (\G,\I);\draw[dashed,line width=0.1mm] (\I,\A) -- (\I,\I);
            \draw[solid,line width=0.5mm] (\A,\A)--(\C,\C)  -- (\E,\A) -- (\G,\C);
         \end{tikzpicture}~~
  \begin{tikzpicture}[scale=0.15]
            \draw (\A,\A)-- (\I,\A);
           
             \draw[dashed,line width=0.1mm] (\A,\A)-- (\I,\A);
              \draw[dashed,line width=0.1mm] (\A,\C)-- (\I,\C);
               \draw[dashed,line width=0.1mm] (\A,\E)-- (\I,\E);
               \draw[dashed,line width=0.1mm] (\A,\G)-- (\I,\G);
               \draw[dashed,line width=0.1mm] (\A,\I)-- (\I,\I);
            \draw (\A,\A) -- (\A,\I);
             \draw[dashed,line width=0.1mm] (\C,\A) -- (\C,\I);\draw[dashed,line width=0.1mm] (\E,\A) -- (\E,\I);\draw[dashed,line width=0.1mm] (\G,\A) -- (\G,\I);\draw[dashed,line width=0.1mm] (\I,\A) -- (\I,\I);
            \draw[solid,line width=0.5mm] (\A,\A)--(\C,\C)  -- (\E,\E)--(\G,\C);
         \end{tikzpicture}~~
          \begin{tikzpicture}[scale=0.15]
            \draw (\A,\A)-- (\I,\A);
            
             \draw[dashed,line width=0.1mm] (\A,\A)-- (\I,\A);
              \draw[dashed,line width=0.1mm] (\A,\C)-- (\I,\C);
               \draw[dashed,line width=0.1mm] (\A,\E)-- (\I,\E);
               \draw[dashed,line width=0.1mm] (\A,\G)-- (\I,\G);
               \draw[dashed,line width=0.1mm] (\A,\I)-- (\I,\I);
            \draw (\A,\A) -- (\A,\I);
             \draw[dashed,line width=0.1mm] (\C,\A) -- (\C,\I);\draw[dashed,line width=0.1mm] (\E,\A) -- (\E,\I);\draw[dashed,line width=0.1mm] (\G,\A) -- (\G,\I);\draw[dashed,line width=0.1mm] (\I,\A) -- (\I,\I);
            \draw[solid,line width=0.5mm] (\A,\A)--(\C,\C)  -- (\E,\E) -- (\G,\E);
         \end{tikzpicture}~~
          \begin{tikzpicture}[scale=0.15]
            \draw (\A,\A)-- (\I,\A);
           
             \draw[dashed,line width=0.1mm] (\A,\A)-- (\I,\A);
              \draw[dashed,line width=0.1mm] (\A,\C)-- (\I,\C);
               \draw[dashed,line width=0.1mm] (\A,\E)-- (\I,\E);
               \draw[dashed,line width=0.1mm] (\A,\G)-- (\I,\G);
               \draw[dashed,line width=0.1mm] (\A,\I)-- (\I,\I);
            \draw (\A,\A) -- (\A,\I);
             \draw[dashed,line width=0.1mm] (\C,\A) -- (\C,\I);\draw[dashed,line width=0.1mm] (\E,\A) -- (\E,\I);\draw[dashed,line width=0.1mm] (\G,\A) -- (\G,\I);\draw[dashed,line width=0.1mm] (\I,\A) -- (\I,\I);
            \draw[solid,line width=0.5mm] (\A,\A)--(\C,\C)  -- (\E,\E) -- (\G,\G);
         \end{tikzpicture}~~
          \begin{tikzpicture}[scale=0.15]
            \draw (\A,\A)-- (\I,\A);
           
             \draw[dashed,line width=0.1mm] (\A,\A)-- (\I,\A);
              \draw[dashed,line width=0.1mm] (\A,\C)-- (\I,\C);
               \draw[dashed,line width=0.1mm] (\A,\E)-- (\I,\E);
               \draw[dashed,line width=0.1mm] (\A,\G)-- (\I,\G);
               \draw[dashed,line width=0.1mm] (\A,\I)-- (\I,\I);
            \draw (\A,\A) -- (\A,\I);
             \draw[dashed,line width=0.1mm] (\C,\A) -- (\C,\I);\draw[dashed,line width=0.1mm] (\E,\A) -- (\E,\I);\draw[dashed,line width=0.1mm] (\G,\A) -- (\G,\I);\draw[dashed,line width=0.1mm] (\I,\A) -- (\I,\I);
            \draw[solid,line width=0.5mm] (\A,\A)--(\C,\C)  --(\E,\E)-- (\G,\A);
         \end{tikzpicture}~~
         \begin{tikzpicture}[scale=0.15]
            \draw (\A,\A)-- (\I,\A);
            
             \draw[dashed,line width=0.1mm] (\A,\A)-- (\I,\A);
              \draw[dashed,line width=0.1mm] (\A,\C)-- (\I,\C);
               \draw[dashed,line width=0.1mm] (\A,\E)-- (\I,\E);
               \draw[dashed,line width=0.1mm] (\A,\G)-- (\I,\G);
               \draw[dashed,line width=0.1mm] (\A,\I)-- (\I,\I);
            \draw (\A,\A) -- (\A,\I);
             \draw[dashed,line width=0.1mm] (\C,\A) -- (\C,\I);\draw[dashed,line width=0.1mm] (\E,\A) -- (\E,\I);\draw[dashed,line width=0.1mm] (\G,\A) -- (\G,\I);\draw[dashed,line width=0.1mm] (\I,\A) -- (\I,\I);
            \draw[solid,line width=0.5mm] (\A,\A)--(\C,\C)--(\E,\C)--(\G,\E) ;
         \end{tikzpicture}~~
         \begin{tikzpicture}[scale=0.15]
            \draw (\A,\A)-- (\I,\A);
           
             \draw[dashed,line width=0.1mm] (\A,\A)-- (\I,\A);
              \draw[dashed,line width=0.1mm] (\A,\C)-- (\I,\C);
               \draw[dashed,line width=0.1mm] (\A,\E)-- (\I,\E);
               \draw[dashed,line width=0.1mm] (\A,\G)-- (\I,\G);
               \draw[dashed,line width=0.1mm] (\A,\I)-- (\I,\I);
            \draw (\A,\A) -- (\A,\I);
             \draw[dashed,line width=0.1mm] (\C,\A) -- (\C,\I);\draw[dashed,line width=0.1mm] (\E,\A) -- (\E,\I);\draw[dashed,line width=0.1mm] (\G,\A) -- (\G,\I);\draw[dashed,line width=0.1mm] (\I,\A) -- (\I,\I);
            \draw[solid,line width=0.5mm] (\A,\A)--(\C,\A)  -- (\E,\C)--(\G,\A) ;
         \end{tikzpicture}~~
         \begin{tikzpicture}[scale=0.15]
            \draw (\A,\A)-- (\I,\A);
           
             \draw[dashed,line width=0.1mm] (\A,\A)-- (\I,\A);
              \draw[dashed,line width=0.1mm] (\A,\C)-- (\I,\C);
               \draw[dashed,line width=0.1mm] (\A,\E)-- (\I,\E);
               \draw[dashed,line width=0.1mm] (\A,\G)-- (\I,\G);
               \draw[dashed,line width=0.1mm] (\A,\I)-- (\I,\I);
            \draw (\A,\A) -- (\A,\I);
             \draw[dashed,line width=0.1mm] (\C,\A) -- (\C,\I);\draw[dashed,line width=0.1mm] (\E,\A) -- (\E,\I);\draw[dashed,line width=0.1mm] (\G,\A) -- (\G,\I);\draw[dashed,line width=0.1mm] (\I,\A) -- (\I,\I);
            \draw[solid,line width=0.5mm] (\A,\A)--(\C,\A)  -- (\E,\C) -- (\G,\C);
         \end{tikzpicture}~~
         \begin{tikzpicture}[scale=0.15]
            \draw (\A,\A)-- (\I,\A);
           
             \draw[dashed,line width=0.1mm] (\A,\A)-- (\I,\A);
              \draw[dashed,line width=0.1mm] (\A,\C)-- (\I,\C);
               \draw[dashed,line width=0.1mm] (\A,\E)-- (\I,\E);
               \draw[dashed,line width=0.1mm] (\A,\G)-- (\I,\G);
               \draw[dashed,line width=0.1mm] (\A,\I)-- (\I,\I);
            \draw (\A,\A) -- (\A,\I);
             \draw[dashed,line width=0.1mm] (\C,\A) -- (\C,\I);\draw[dashed,line width=0.1mm] (\E,\A) -- (\E,\I);\draw[dashed,line width=0.1mm] (\G,\A) -- (\G,\I);\draw[dashed,line width=0.1mm] (\I,\A) -- (\I,\I);
            \draw[solid,line width=0.5mm] (\A,\A)--(\C,\A)  -- (\E,\C) -- (\G,\E);
         \end{tikzpicture}
               \end{center}
         \caption{ The 9 PMAP  in $\mathcal{M}_2$. Notice that two paths end on the $x$-axis, three paths end at height 1, three paths end at height 2, and one path end at height 3.}
         \label{figmot}
\end{figure}
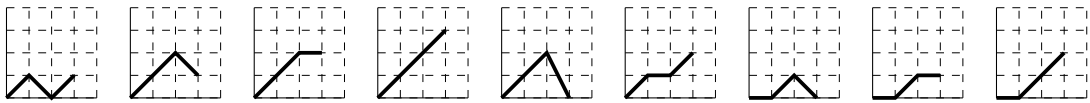

The first few terms of the series expansion of $\mathit{Total}(z,1)$ are  $1 + 2z + 4z^2 + 9z^3 + 21z^4 + 51z^5 + 127z^6 + 323z^7 + 835z^8 + 2188z^9+ O(z^{10})$, which correspond to a shift of the sequence \oeis{A001006} in \cite{oeis} that counts the Motzkin paths of a given length. See Figure \ref{figmot} for an illustration of the 9 paths of length 3.

\begin{cor} The g.f. that counts the MAP with respect to the length is given by 
$$\mathit{Total}(z,0)=\frac{1}{zr}.$$
\end{cor}

The first few terms of the series expansion of $\mathit{Total}(z,0)$ are
$1 + z + z^2 + 2z^3 + 4z^4 + 9z^5 + 21z^6 + 51z^7 + 127z^8 + 323z^9 + O(x^{10})$
which correspond to a shift of the sequence \oeis{A001006} in \cite{oeis} that counts Motzkin paths of a given length.

\subsection{PMAP in $\mathcal{M}'_2$ - From right to left}
\begin{figure}[h]
 \begin{center}
        \begin{tikzpicture}[scale=0.15]
            \draw (\A,\A)-- (38,\A);
             \draw[dashed,line width=0.1mm] (\A,\E)-- (\Ze,\E);
              \draw[dashed,line width=0.1mm] (\A,\C)-- (\Ze,\C);
               \draw[dashed,line width=0.1mm] (\A,\G)-- (\Ze,\G);
               \draw[dashed,line width=0.1mm] (\A,\I)-- (\Ze,\I);
              \draw[dashed,line width=0.1mm] (\A,\K)-- (\Ze,\K);
               \draw[dashed,line width=0.1mm] (\A,\M)-- (\Ze,\M);
            \draw (\A,\A) -- (\A,\O);
             \draw[dashed,line width=0.1mm] (\C,\A) -- (\C,\M);\draw[dashed,line width=0.1mm] (\E,\A) -- (\E,\M);\draw[dashed,line width=0.1mm] (\G,\A) -- (\G,\M);
             \draw[dashed,line width=0.1mm] (\I,\A) -- (\I,\M);\draw[dashed,line width=0.1mm] (\K,\A) -- (\K,\M);\draw[dashed,line width=0.1mm] (\M,\A) -- (\M,\M);
             \draw[dashed,line width=0.1mm] (\O,\A) -- (\O,\M);\draw[dashed,line width=0.1mm] (\Q,\A) -- (\Q,\M);\draw[dashed,line width=0.1mm] (\S,\A) -- (\S,\M);
             \draw[dashed,line width=0.1mm] (\U,\A) -- (\U,\M);\draw[dashed,line width=0.1mm] (\W,\A) -- (\W,\M);\draw[dashed,line width=0.1mm] (\Y,\A) -- (\Y,\M);
             \draw[dashed,line width=0.1mm] (\ZZ,\A) -- (\ZZ,\M);
             \draw[dashed,line width=0.1mm] (\Za,\A) -- (\Za,\M);
             \draw[dashed,line width=0.1mm] (\Zb,\A) -- (\Zb,\M);
             \draw[dashed,line width=0.1mm] (\Zc,\A) -- (\Zc,\M);
             \draw[dashed,line width=0.1mm] (\Zd,\A) -- (\Zd,\M);
             \draw[dashed,line width=0.1mm] (\Ze,\A) -- (\Ze,\M);
            \draw[solid,line width=0.4mm] (\A,\A)--(\C,\G)  -- (\E,\E) -- (\G,\C) --(\I,\I)-- (\K,\G) -- (\M,\E) -- (\O,\E) -- (\Q,\C)  -- (\S,\A)--(\U,\G)--(\W,\E)--(\Y,\E)--(\Za,\A) -- (\Zb,\G) --(\Ze,\A);

         \end{tikzpicture}\qquad 
 \begin{tikzpicture}[scale=0.15]
            \draw (\A,\A)-- (38,\A);
             \draw[dashed,line width=0.1mm] (\A,\E)-- (\Ze,\E);
              \draw[dashed,line width=0.1mm] (\A,\C)-- (\Ze,\C);
               \draw[dashed,line width=0.1mm] (\A,\G)-- (\Ze,\G);
               \draw[dashed,line width=0.1mm] (\A,\I)-- (\Ze,\I);
              \draw[dashed,line width=0.1mm] (\A,\K)-- (\Ze,\K);
               \draw[dashed,line width=0.1mm] (\A,\M)-- (\Ze,\M);
            \draw (\A,\A) -- (\A,\O);
             \draw[dashed,line width=0.1mm] (\C,\A) -- (\C,\M);\draw[dashed,line width=0.1mm] (\E,\A) -- (\E,\M);\draw[dashed,line width=0.1mm] (\G,\A) -- (\G,\M);
             \draw[dashed,line width=0.1mm] (\I,\A) -- (\I,\M);\draw[dashed,line width=0.1mm] (\K,\A) -- (\K,\M);\draw[dashed,line width=0.1mm] (\M,\A) -- (\M,\M);
             \draw[dashed,line width=0.1mm] (\O,\A) -- (\O,\M);\draw[dashed,line width=0.1mm] (\Q,\A) -- (\Q,\M);\draw[dashed,line width=0.1mm] (\S,\A) -- (\S,\M);
             \draw[dashed,line width=0.1mm] (\U,\A) -- (\U,\M);\draw[dashed,line width=0.1mm] (\W,\A) -- (\W,\M);\draw[dashed,line width=0.1mm] (\Y,\A) -- (\Y,\M);
             \draw[dashed,line width=0.1mm] (\ZZ,\A) -- (\ZZ,\M);
             \draw[dashed,line width=0.1mm] (\Za,\A) -- (\Za,\M);
             \draw[dashed,line width=0.1mm] (\Zb,\A) -- (\Zb,\M);
             \draw[dashed,line width=0.1mm] (\Zc,\A) -- (\Zc,\M);
             \draw[dashed,line width=0.1mm] (\Zd,\A) -- (\Zd,\M);
             \draw[dashed,line width=0.1mm] (\Ze,\A) -- (\Ze,\M);
            \draw[solid,line width=0.4mm](\A,\A)--(\C,\G)  -- (\E,\E) -- (\G,\C) --(\I,\I)-- (\K,\G) -- (\M,\E) -- (\O,\E) -- (\Q,\C)  -- (\S,\A)--(\U,\G)--(\W,\E)--(\Y,\E)--(\Za,\A) -- (\Zb,\G) --(\Zc,\E)--(\Zd,\G)--(\Ze,\E);;

         \end{tikzpicture}
               \end{center}
         \caption{  The left drawing shows a Motzkin path with air pockets of length $18$ read from right to left. The right drawing shows a partial Motzkin path with air pockets of  length $18$  ending at height $2$ and reaf from right to left.}
         \label{fig4}
\end{figure}
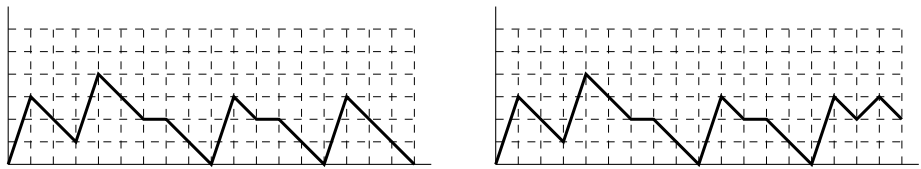

Here, we consider the paths of the previous section, but we read them from right to left. This means that down steps become up steps  and {\it vice versa}, and horizontal steps are unchanged, which implies that any up step or horizontal step (except the first step of the path) is preceded by a down step. See Figure \ref{fig4} for two examples of such paths.

 So, we easily obtain the following equations:
\begin{equation}\left\{\begin{array}{l}
f_0=1,\mbox{ and } f_k=z(1+g_0+g_1+\ldots +g_{k-1}), \quad k\geq 1,\\
g_k=zf_{k+1}+zg_{k+1}+zh_{k+1},\quad k\geq 0,\\
h_0=z+zg_0,\mbox{ and } h_k=zg_k,\quad k\geq 1.\\
\end{array}\right.\label{equ4}
\end{equation}

Using the same notations as in the previous sections, and summing the recursions in (\ref{equ4}), we have:
\begin{align*}
F(u)&=1+z\sum\limits_{k\geq 1}u^kf_{k-1}=1+z\sum\limits_{k\geq 1}(1+g_0+\ldots +g_{k-1})u^k\\
&=1+\frac{zu}{1-u}(1+G(u)),\\
G(u)&=z\sum\limits_{k\geq 0}(f_{k+1}+g_{k+1}+h_{k+1})u^k\\
&=\frac{z}{u}(F(u)-F(0)+G(u)-G(0)+H(u)-H(0)),\\
H(u)&=z+zG(u).
\end{align*}

Notice that $F(0)=1$ and $H(0)=z+zG(0)$ by the third relation. Now, setting $g_0:=G(0)$  and solving these functional equations, we deduce
$$F(u)={\frac {{\it g_0}\,u{z}^{3}+{\it g_0}\,u{z}^{2}+u{z}^{3}-{u}^{2}z+u{z}^{
2}+{u}^{2}-zu+{z}^{2}-u+z}{{u}^{2}-zu+{z}^{2}-u+z}}
,$$
$$G(u)=-{\frac {z \left( {\it g_0}\,uz+{\it g_0}\,u-{\it g_0}\,z+zu-{\it g_0}
 \right) }{{u}^{2}-zu+{z}^{2}-u+z}},$$
$$H(u)=-{\frac {z \left( {\it g_0}\,u{z}^{2}+{\it g_0}\,uz-{\it g_0}\,{z}^{2}+u{
z}^{2}-{\it g_0}\,z-{u}^{2}+zu-{z}^{2}+u-z \right) }{{u}^{2}-zu+{z}^{2}
-u+z}}
.$$

In order to compute $g_0$, we use the kernel method (see~\cite{ban, pro}) on $F(u)$. We can write the denominator (which is a polynomial in $u$ of degree 2), as $(u-r)(u-s)$ with 
$$r=\frac{z+1+\sqrt {-3\,{z}^{2}-2\,z+1}}{2}
,\mbox{ and }
s=\frac{z+1-\sqrt {-3\,{z}^{2}-2\,z+1}}{2}.$$
Plugging $u=s$ (which have a Taylor expansion at $z=0$) in $F(u)(u-r)(u-s)$, we obtain the equation 
$${\it g_0}\,s{z}^{3}+{\it g_0}\,s{z}^{2}+s{z}^{3}-{s}^{2}z+s{z}^{
2}+{s}^{2}-zs+{z}^{2}-s+z=0.$$
Using $rs=z(1+z)$, we deduce 
 $$g_0=\frac{1-r}{r}
.$$
Finally, after simplifying by the factor $(u-s)$ in the numerators and denominators, we obtain 
$$F(u)=\frac{u(z-1)+r}{r-u}
,\quad G(u)=\frac{1-r}{r-u},\quad \mbox{ and }\quad  H(u)=\frac{z(1-u)}{r-u}
,$$
which implies that
\begin{align}f_k&=[u^k] F(u)= \frac{1}{r^k}+(1-[k=0])\cdot\frac{z-1}{r^k},\\
g_k&=[u^k] G(u)= \frac{1-r}{r^{k+1}},\\
h_k&=[u^k] H(u)=\frac{z}{r^{k+1}}-(1-[k=0])\cdot\frac{z}{r^{k}}.
\end{align}

\begin{thm} The bivariate generating function for the total number of PMAP with respect to the length and the height of the end-point is given by
$$\mathit{Total}(z,u)={\frac {1-u+z}{r-u}}
.$$

$$[u^k]total(z,u)=\frac{z+1}{r^{k+1}}-(1-[k=0])\frac{1}{r^{k}}.$$

Finally, setting $t(n,k)=[z^n][u^k]\mathit{Total}(z,u)$, we have for $n\geq 1$, $k\geq 1$,  $$t(n,k)=t(n,k-1)-t(n-1,k)+t(n-2,k+1)+t(n-1,k+1)
,$$
and setting $t_n:=t(n,0)$,  we have $$t_n=t_{n-1}+\sum\limits_{k=1}^{n-2}t_kt_{n-1-k}.$$
\end{thm}
\noindent {\it Proof.} The proof are obtained {\it mutatis mutandis} as for the previous theorems.
\hfill $\Box$

Let $\mathcal{T}$ be the infinite matrix $\mathcal{T}:=[t(n,k)]_{n\geq 0, k\geq 0}$. The first few rows of the matrix $\mathcal{T}$ are
\[\mathcal{T}=
\left(
\begin{array}{cccccccc}
 1 & 0 & 0 & 0 & 0 & 0 & 0 & 0  \\
 1 & 1 & 1 & 1 & 1 & 1 & 1 & 1  \\
 1 & 1 & 1 & 1 & 1 & 1 & 1 & 1  \\
 2 & 3 & 4 & 5 & 6 & 7 & 8 & 9  \\
 4 & 6 & 8 & 10 & 12 & 14 & 16 & 18  \\
 9 & 15 & 22 & 30 & 39 & 49 & 60 & 72 \\
 21 & 36 & 54 & 75 & 99 & 126 & 156 & 189  \\
 51 & 91 & 142 & 205 & 281 & 371 & 476 & 597  \\
 127 & 232 & 370 & 545 & 761 & 1022 & 1332 & 1695  \\
\end{array}
\right).\]

Since there is an infinite number of PMAP of length $n$, we do not provide an ordinary generating function (with
respect to the length) for these paths. So, we get around this  by counting  PMAP ending on a point $(x,n-x)$ for a given $n\geq 0$.

\begin{cor} The g.f. that counts the partial PMAP ending on the line $y=n-x$ is given by 
$$\mathit{Total}(z,z)=\frac{1}{r-z}.$$
\end{cor}

\begin{figure}[h]
 \begin{center}
        \begin{tikzpicture}[scale=0.15]
            \draw (\A,\A)-- (\I,\A);
            \draw[red] (\A,\I)--(\I,\A);
             \draw[dashed,line width=0.1mm] (\A,\A)-- (\I,\A);
              \draw[dashed,line width=0.1mm] (\A,\C)-- (\I,\C);
               \draw[dashed,line width=0.1mm] (\A,\E)-- (\I,\E);
               \draw[dashed,line width=0.1mm] (\A,\G)-- (\I,\G);
               \draw[dashed,line width=0.1mm] (\A,\I)-- (\I,\I);
            \draw (\A,\A) -- (\A,\I);
             \draw[dashed,line width=0.1mm] (\C,\A) -- (\C,\I);\draw[dashed,line width=0.1mm] (\E,\A) -- (\E,\I);\draw[dashed,line width=0.1mm] (\G,\A) -- (\G,\I);\draw[dashed,line width=0.1mm] (\I,\A) -- (\I,\I);
            \draw[solid,line width=0.5mm] (\A,\A)--(\C,\C)  -- (\E,\A) -- (\G,\C) --(\I,\A);
         \end{tikzpicture}~~
  \begin{tikzpicture}[scale=0.15]
            \draw (\A,\A)-- (\I,\A);
            \draw[red] (\A,\I)--(\I,\A);
             \draw[dashed,line width=0.1mm] (\A,\A)-- (\I,\A);
              \draw[dashed,line width=0.1mm] (\A,\C)-- (\I,\C);
               \draw[dashed,line width=0.1mm] (\A,\E)-- (\I,\E);
               \draw[dashed,line width=0.1mm] (\A,\G)-- (\I,\G);
               \draw[dashed,line width=0.1mm] (\A,\I)-- (\I,\I);
            \draw (\A,\A) -- (\A,\I);
             \draw[dashed,line width=0.1mm] (\C,\A) -- (\C,\I);\draw[dashed,line width=0.1mm] (\E,\A) -- (\E,\I);\draw[dashed,line width=0.1mm] (\G,\A) -- (\G,\I);\draw[dashed,line width=0.1mm] (\I,\A) -- (\I,\I);
            \draw[solid,line width=0.5mm] (\A,\A)--(\C,\C)  -- (\E,\A) -- (\G,\C);
         \end{tikzpicture}~~
          \begin{tikzpicture}[scale=0.15]
            \draw (\A,\A)-- (\I,\A);
            \draw[red] (\A,\I)--(\I,\A);
             \draw[dashed,line width=0.1mm] (\A,\A)-- (\I,\A);
              \draw[dashed,line width=0.1mm] (\A,\C)-- (\I,\C);
               \draw[dashed,line width=0.1mm] (\A,\E)-- (\I,\E);
               \draw[dashed,line width=0.1mm] (\A,\G)-- (\I,\G);
               \draw[dashed,line width=0.1mm] (\A,\I)-- (\I,\I);
            \draw (\A,\A) -- (\A,\I);
             \draw[dashed,line width=0.1mm] (\C,\A) -- (\C,\I);\draw[dashed,line width=0.1mm] (\E,\A) -- (\E,\I);\draw[dashed,line width=0.1mm] (\G,\A) -- (\G,\I);\draw[dashed,line width=0.1mm] (\I,\A) -- (\I,\I);
            \draw[solid,line width=0.5mm] (\A,\A)--(\C,\E)  -- (\E,\C) -- (\G,\C)--(\I,\A);
         \end{tikzpicture}~~
          \begin{tikzpicture}[scale=0.15]
            \draw (\A,\A)-- (\I,\A);
            \draw[red] (\A,\I)--(\I,\A);
             \draw[dashed,line width=0.1mm] (\A,\A)-- (\I,\A);
              \draw[dashed,line width=0.1mm] (\A,\C)-- (\I,\C);
               \draw[dashed,line width=0.1mm] (\A,\E)-- (\I,\E);
               \draw[dashed,line width=0.1mm] (\A,\G)-- (\I,\G);
               \draw[dashed,line width=0.1mm] (\A,\I)-- (\I,\I);
            \draw (\A,\A) -- (\A,\I);
             \draw[dashed,line width=0.1mm] (\C,\A) -- (\C,\I);\draw[dashed,line width=0.1mm] (\E,\A) -- (\E,\I);\draw[dashed,line width=0.1mm] (\G,\A) -- (\G,\I);\draw[dashed,line width=0.1mm] (\I,\A) -- (\I,\I);
            \draw[solid,line width=0.5mm] (\A,\A)--(\C,\E)  -- (\E,\C) -- (\G,\A)--(\I,\A);
         \end{tikzpicture}~~
          \begin{tikzpicture}[scale=0.15]
            \draw (\A,\A)-- (\I,\A);
            \draw[red] (\A,\I)--(\I,\A);
             \draw[dashed,line width=0.1mm] (\A,\A)-- (\I,\A);
              \draw[dashed,line width=0.1mm] (\A,\C)-- (\I,\C);
               \draw[dashed,line width=0.1mm] (\A,\E)-- (\I,\E);
               \draw[dashed,line width=0.1mm] (\A,\G)-- (\I,\G);
               \draw[dashed,line width=0.1mm] (\A,\I)-- (\I,\I);
            \draw (\A,\A) -- (\A,\I);
             \draw[dashed,line width=0.1mm] (\C,\A) -- (\C,\I);\draw[dashed,line width=0.1mm] (\E,\A) -- (\E,\I);\draw[dashed,line width=0.1mm] (\G,\A) -- (\G,\I);\draw[dashed,line width=0.1mm] (\I,\A) -- (\I,\I);
            \draw[solid,line width=0.5mm] (\A,\A)--(\C,\E)  -- (\E,\C) -- (\G,\C);
         \end{tikzpicture}~~
         \begin{tikzpicture}[scale=0.15]
            \draw (\A,\A)-- (\I,\A);
            \draw[red] (\A,\I)--(\I,\A);
             \draw[dashed,line width=0.1mm] (\A,\A)-- (\I,\A);
              \draw[dashed,line width=0.1mm] (\A,\C)-- (\I,\C);
               \draw[dashed,line width=0.1mm] (\A,\E)-- (\I,\E);
               \draw[dashed,line width=0.1mm] (\A,\G)-- (\I,\G);
               \draw[dashed,line width=0.1mm] (\A,\I)-- (\I,\I);
            \draw (\A,\A) -- (\A,\I);
             \draw[dashed,line width=0.1mm] (\C,\A) -- (\C,\I);\draw[dashed,line width=0.1mm] (\E,\A) -- (\E,\I);\draw[dashed,line width=0.1mm] (\G,\A) -- (\G,\I);\draw[dashed,line width=0.1mm] (\I,\A) -- (\I,\I);
            \draw[solid,line width=0.5mm] (\A,\A)--(\C,\G) ;
         \end{tikzpicture}~~
         \begin{tikzpicture}[scale=0.15]
            \draw (\A,\A)-- (\I,\A);
            \draw[red] (\A,\I)--(\I,\A);
             \draw[dashed,line width=0.1mm] (\A,\A)-- (\I,\A);
              \draw[dashed,line width=0.1mm] (\A,\C)-- (\I,\C);
               \draw[dashed,line width=0.1mm] (\A,\E)-- (\I,\E);
               \draw[dashed,line width=0.1mm] (\A,\G)-- (\I,\G);
               \draw[dashed,line width=0.1mm] (\A,\I)-- (\I,\I);
            \draw (\A,\A) -- (\A,\I);
             \draw[dashed,line width=0.1mm] (\C,\A) -- (\C,\I);\draw[dashed,line width=0.1mm] (\E,\A) -- (\E,\I);\draw[dashed,line width=0.1mm] (\G,\A) -- (\G,\I);\draw[dashed,line width=0.1mm] (\I,\A) -- (\I,\I);
            \draw[solid,line width=0.5mm] (\A,\A)--(\C,\G)  -- (\E,\E) ;
         \end{tikzpicture}~~
         \begin{tikzpicture}[scale=0.15]
            \draw (\A,\A)-- (\I,\A);
            \draw[red] (\A,\I)--(\I,\A);
             \draw[dashed,line width=0.1mm] (\A,\A)-- (\I,\A);
              \draw[dashed,line width=0.1mm] (\A,\C)-- (\I,\C);
               \draw[dashed,line width=0.1mm] (\A,\E)-- (\I,\E);
               \draw[dashed,line width=0.1mm] (\A,\G)-- (\I,\G);
               \draw[dashed,line width=0.1mm] (\A,\I)-- (\I,\I);
            \draw (\A,\A) -- (\A,\I);
             \draw[dashed,line width=0.1mm] (\C,\A) -- (\C,\I);\draw[dashed,line width=0.1mm] (\E,\A) -- (\E,\I);\draw[dashed,line width=0.1mm] (\G,\A) -- (\G,\I);\draw[dashed,line width=0.1mm] (\I,\A) -- (\I,\I);
            \draw[solid,line width=0.5mm] (\A,\A)--(\C,\G)  -- (\E,\E) -- (\G,\C);
         \end{tikzpicture}~~
         \begin{tikzpicture}[scale=0.15]
            \draw (\A,\A)-- (\I,\A);
            \draw[red] (\A,\I)--(\I,\A);
             \draw[dashed,line width=0.1mm] (\A,\A)-- (\I,\A);
              \draw[dashed,line width=0.1mm] (\A,\C)-- (\I,\C);
               \draw[dashed,line width=0.1mm] (\A,\E)-- (\I,\E);
               \draw[dashed,line width=0.1mm] (\A,\G)-- (\I,\G);
               \draw[dashed,line width=0.1mm] (\A,\I)-- (\I,\I);
            \draw (\A,\A) -- (\A,\I);
             \draw[dashed,line width=0.1mm] (\C,\A) -- (\C,\I);\draw[dashed,line width=0.1mm] (\E,\A) -- (\E,\I);\draw[dashed,line width=0.1mm] (\G,\A) -- (\G,\I);\draw[dashed,line width=0.1mm] (\I,\A) -- (\I,\I);
            \draw[solid,line width=0.5mm] (\A,\A)--(\C,\G)  -- (\E,\E) -- (\G,\C)--(\I,\A);
         \end{tikzpicture}
               \end{center}
         \caption{ The 9 PMAP ending on the line $y=4-x$ in $\mathcal{M}'_2$. Notice that four paths end on the $x$-axis.}
         \label{figlast}
\end{figure}
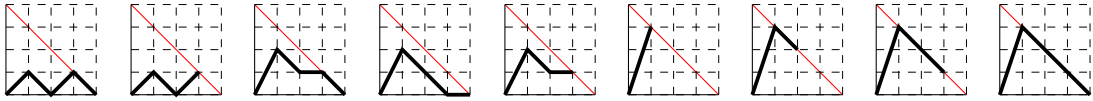

The first few terms of the series expansion of $\mathit{Total}(z,z)$ are  $1 + z + 2z^2 + 4z^3 + 9z^4 + 21z^5 + 51z^6 + 127z^7 + 323z^8 + 835z^9+ O(z^{10})$, which correspond to the sequence \oeis{A001006} in \cite{oeis} that counts the  Motzkin paths with respect to the length. See Figure \ref{figlast} for the illustration of the 9 PMAP in $\mathcal{M}'_2$ ending on the line $y=4-x$.

Notice that we obviously retrieve the results of Corollary 5, i.e., the g.f. $\mathit{Total}(z,0)$ that counts the MAP with respect to the length is also a shift of the Motzkin sequence \oeis{A001006} in \cite{oeis}.

\section{A Riordan array point of view}
In this section, we make links between the previous matrices $\mathcal{T}=[t_{n,k}]_{n\geq 0, k \geq 0}$ and some Riordan arrays or `almost' Riordan arrays. We first give a short background on Riordan arrays  \cite{Barry,Bur,Riordan}.

An infinite column vector $(a_0,a_1,\dots)^T$ has generating function $f(z)$ if $f(z)=\sum_{n\ge0} a_nz^n$. A {\em Riordan array} is an infinite lower triangular matrix whose $k$-th column has generating function
$g(z)f(z)^k$ for all $k \ge 0$, for some formal power series  $g(z)$ and $f(z)$, with $g(0) \neq 0$, $f(0)=0$, and $f'(0)\neq 0$. Such a Riordan array is denoted by  $(g(z),f(z))$.
If we multiply this matrix by a column vector  $(c_0, c_1, \dots)^T$ having generating function $h(z)$, then the resulting column
vector has generating function $g(z)h(f(z))$. This property is known as the {\it fundamental theorem of Riordan arrays}.

The product of two Riordan arrays $(g(z),f(z))$ and $(h(z),l(z))$ is defined by
\begin{equation}\label{eq:product} (g(z),f(z))*(h(z),l(z))=\left(g(z)h(f(z)), l(f(z))\right). \end{equation}
Under the operation $``*"$, the set of all Riordan arrays is a group  \cite{Riordan}. The identity element is $I=(1,z)$, and the inverse of $(g(z),f(z))$ is
 \begin{equation}\label{invRiordan}
 (g(z), f(z))^{-1}=\left(1/\left(g\circ \overline{f}\right)(z), \overline{f}(z)\right),
\end{equation}
  where $\overline{f}(z)$ denotes the compositional inverse of $f(z)$.
  
 Finally, we will say that a matrix $\mathcal{M}$ is the {\it rectification of the Riordan array} $(g(z),f(z))$ whenever the bivariate generating function of  $\mathcal{M}$ equals $$\frac{g(z)}{1-u\frac{f(z)}{z}}.$$
  
\subsection{Comment on Section 2.1}
\begin{prop}
The matrix $\mathcal{T}=[t(n,k)]_{n\geq 0, k \geq 0}$ is a Riordan array defined by 
$$\left(\frac{1}{1-z^2}C\left(\frac{z(1-z+z^2)}{(1-z^2)^2}\right),\frac{z}{1-z^2}C\left(\frac{z(1-z+z^2)}{(1-z^2)^2}\right)\right)$$
where $C(z)=\frac{1-\sqrt{1-4z}}{2z}$ is the generating function of the Catalan numbers
$c_n=\frac{1}{n+1}\binom{2n}{n}$.
\end{prop}
\noindent {\it Proof.}
 Indeed, since $1+rs-z=\frac{1}{z}$, we have 
\begin{align*}
[u^k]\mathit{Total}(z,u)&=\frac{1+rs-z}{r^{k+1}}=\frac{1}{zr^{k+1}}=\left(\frac{1}{zr}\right)\left(\frac{1}{r}\right)^k\\
&=\frac{2}{1-z^2+\sqrt{1-4z+2z^2-4z^3+z^4}}\left(\frac{2z}{1-z^2+\sqrt{1-4z+2z^2-4z^3+z^4}}\right)^k.
\end{align*}
Therefore, the array $\mathcal{T}$ satisfies
\begin{align*}\mathcal{T}&=\left(\frac{2}{1-z^2+\sqrt{1-4z+2z^2-4z^3+z^4}},
\frac{2z}{1-z^2+\sqrt{1-4z+2z^2-4z^3+z^4}}\right)\\
&=\left(\frac{1-z^2-\sqrt{1-4z+2z^2-4z^3+z^4}}{2(1-z+z^2)},
\frac{z(1-z^2-\sqrt{1-4z+2z^2-4z^3+z^4})}{2(1-z+z^2)}\right)\\
&=\left(\frac{1}{1-z^2}C\left(\frac{z(1-z+z^2)}{(1-z^2)^2}\right),\frac{z}{1-z^2}C\left(\frac{z(1-z+z^2)}{(1-z^2)^2}\right)\right),\end{align*}
where $C(z)=\frac{1-\sqrt{1-4z}}{2z}$ is the generating function of the Catalan numbers
$c_n=\frac{1}{n+1}\binom{2n}{n}$.
\hfill $\Box$
\medskip

As a consequence, we have
$$t(n,k)=\sum_{i=0}^n \binom{k+\frac{i}{2}}{\frac{i}{2}}\frac{1+(-1)^i}{2}\sum_{j=0}^{n-i}A_{j,k}B_{n-i-j,k},$$
where
\begin{align*}A_{n,k}&=\sum_{j=0}^n C_{j,k} D_{n,j},\\
B_{n,k}&=\sum_{j=0}^{2k} \binom{2k}{j}(-1)^j \sum_{i=0}^{n-2j}\binom{k+i-1}{i}\binom{i}{n-2j-i}(-1)^{n-i},\\
C_{n,k}&=\frac{k+1}{n+1}\binom{2n-k}{n-k},\\
D_{n,k}&=\sum_{j=0}^k \binom{k}{j}(-1)^j \sum_{i=0}^j \binom{j}{i}\binom{2k-1+\frac{n-k-j-i}{2}}{2k-1}(-1)^i
\frac{1+(-1)^{n-k-j-i}}{2}. \end{align*}
The matrix $[C_{n,k}]_{n,k\geq 0}$ is the Riordan array  $(C(z), zC(z))$ (see the Catalan matrix \oeis{A033184} in \cite{oeis}).
The matrix $[D_{n,k}]_{n,k\geq 0}$ is the Riordan array $\left(1,\frac{z(1-z+z^2)}{(1-z^2)^2}\right)$.

The sequence $u_n = \sum_{k=0}^n D_{n,k}c_k$, with generating function $C\left(\frac{z(1-z+z^2)}{(1-z^2)^2}\right)$,  begins 
$$1,1,1,4,11,31,92,281,877,2788,8999,29415,\ldots.$$
Then the sequence $t(n,0)$ (\oeis{A114465}) which begins
$$1, 1, 2, 5, 13, 36, 105, 317, 982, 3105, 9981, 32520, \ldots$$ is the convolution of $(u_n)_{n\geq 0}=1,1,1,4,11,31,\ldots$ and $1,0,1,0,1,0,\ldots$, which means that
$$t(n,0) = \sum_{k=0}^n u_{n-k} \frac{1+(-1)^k}{2}.$$ 

\begin{prop} The general term $t(n,k)$ equals 
$$\scalemath{0.75}{\sum\limits_{i=0}^{n-k}\binom{k+\frac{n-k-i}{2}}{\frac{n-k-i}{2}}\frac{1+(-1)^{n-k-i}}{2}\sum\limits_{j=0}^iC_{k+j,k}\sum\limits_{m=0}^{j}\binom{j}{m} (-1)^m\sum\limits_{p=0}^m\binom{m}{p}(-1)^p\binom{2j-1+\frac{i-j-m-p}{2}}{\frac{i-j-m-p}{2}}\frac{1+(-1)^{i-j-m-p}}{2}},$$
where $C_{n,k}=\frac{k+1}{n+1}\binom{2n-k}{n-k}$ is the general term of the Catalan matrix (\oeis{A033184}).
\end{prop}
\noindent {\it Proof.} Setting  $Z=\frac{z(1-z+z^2)}{(1-z^2)^2}$, we have
\begin{align*}
t(n,k)&=[z^n]z^k\frac{1}{(1-z^2)^{k+1}}C(Z)^{k+1}\\
&=[z^{n-k}]\frac{1}{(1-z^2)^{k+1}}G(Z) \quad \mbox{ where } G(z)=C(z)^{k+1}\\
&=\sum\limits_{i=0}^{n-k}[z^{n-k-i}]\frac{1}{(1-z^2)^{k+1}}[z^i]G(Z)\\
&=\sum_{i=0}^{n-k}[z^{n-k-i}]\frac{1}{(1-z^2)^{k+1}}\sum\limits_{j=0}^i[z^j]C(z)^{k+1}[z^i]Z^j\\
&=\sum\limits_{i=0}^{n-k}[z^{n-k-i}]\frac{1}{(1-z^2)^{k+1}}\sum\limits_{j=0}^{i}[z^j]\frac{1}{z^k}C(z)(zC(z))^k[z^i]Z^j\\
&=\sum\limits_{i=0}^{n-k}\binom{k+\frac{n-k-i}{2}}{\frac{n-k-i}{2}}\frac{1+(-1)^{n-k-i}}{2}\sum\limits_{j=0}^{i}C_{k+j,k}[z^i]Z^j.
\end{align*}
Since we have 
\begin{align*}
[z^i]Z^j&=[z^i]\left(\frac{z(1-z+z^2)}{(1-z^2)^2}\right)^j\\
&=\sum\limits_{m=0}^j\binom{j}{m}(-1)^m\sum\limits_{p=0}^m\binom{m}{p}(-1)^p\binom{2j-1+\frac{i-j-m-p}{2}}{\frac{i-j-m-p}{2}}\frac{1+(-1)^{i-j-m-p}}{2},
\end{align*}
the result follows.
\hfill $\Box$

\subsection{Comment on Section 2.2}

\begin{prop} The matrix $\mathcal{T}=[t(n,k)]_{n\geq 0, k \geq 0}$ can be written  
\[\mathcal{T}=
\left(
\begin{array}{ccccccc}
 1 & 0 & 0 & 0 & 0 & 0 &\cdots \\
 1 & 1 & 1 & 1 & 1 & 1 &\cdots  \\
 2 & 3 & 3 & 3 & 3 & 3 &\cdots  \\
 5 & 8 & 10 & 12 & 14 & 16  &\cdots  \\
 13 & 23 & 33 & 43 & 53 & 63 &\cdots  \\
 36 & 69 & 107 & 149 & 195 & 245 &\cdots \\
 \vdots & \vdots & \vdots & \vdots & \vdots & \vdots &\ddots \\
\end{array}
\right)=A\cdot B\]
where \[ A=\left(
\begin{array}{ccccccc}
 1 & 0 & 0 & 0 & 0 & 0 &\cdots  \\
 1 & 1 & 0 & 0 & 0 & 0 &\cdots \\
 2 & 3 & 0 & 0 & 0 & 0 &\cdots  \\
 5 & 8 & 2 & 0 & 0 & 0  &\cdots  \\
 13 & 23 & 10 & 0 & 0 & 0 &\cdots  \\
 36 & 69 & 38 & 4 & 0 & 0 &\cdots \\
  \vdots & \vdots & \vdots & \vdots & \vdots& \vdots  &\ddots \\
\end{array}
\right)\quad \mbox{and}\quad B=\left(
\begin{array}{ccccccc}
 1 & 0 & 0 & 0 & 0 & 0 &\cdots  \\
 0 & 1 & 1 & 1 & 1 & 1 &\cdots \\
 0 & 0 & 1 & 2 & 3 & 4 &\cdots  \\
 0 & 0 & 0 & 1 & 3 & 6  &\cdots  \\
 0 & 0 & 0 & 0 & 1 & 4 &\cdots  \\
 0 & 0 & 0 & 0 & 0 & 1 &\cdots  \\
  \vdots & \vdots & \vdots & \vdots & \vdots& \vdots  &\ddots \\
\end{array}
\right)
\]
are defined as follows:

 $\bullet$ The matrix $B=[b_{n,k}]_{n,k\geq 0}$ is defined by $b_{0,0}=1$, and $b_{n,0}=b_{0,n}=0$ if $n\geq 1$, and $b_{n,k}=\binom{k-1}{n-1}$ otherwise, which is a kind of Pascal matrix. 

 $\bullet$ The matrix $A=[a_{n,k}]_{n,k\geq 0}$ is the `almost' Riordan array with initial column whose generating function is $g_0(z)$ which is followed by the shifted `stretched' Riordan array 
$$\scalemath{0.85}{\left(\frac{1-3z+z^2-z^3-(1-z)\sqrt{1-4z+2z^2-4z^3+z^4}}{2z^3(1-z+z^2)},\frac{2(1-2z-z^2-\sqrt{1-4z+2z^2-4z^3+z^4})}{4z}\right)} .$$
The sequence $1,3,8,23,69,\ldots$ with generating function $\frac{1-3z+z^2-z^3-(1-z)\sqrt{1-4z+2z^2-4z^3+z^4}}{2z^3(1-z+z^2)}$ is the convolution of $1,1,2,5,13,36,\ldots$ (\oeis{A114465}) and $1,2,4,10,28,\ldots$ (\oeis{A187256}).
\end{prop}

\noindent {\it Proof.} An almost Riordan array of order $1$ is represented by an initial column vector with generating function $g_0(z)$, followed by a vertically shifted Riordan array $(g(z),f(z))$. The bivariate generating function of this matrix is then given by $g_0(z)+zu\frac{g(z)}{1-uf(z)}$. In our case, we have
\begin{align*}
g_0(z)&=\frac{1-z^2-\sqrt{1-4z+2z^2-4z^3+z^4}}{2z(1-z+z^2)},\\
g(z)&=\frac{1-3z+z^2-z^3-(1-z)\sqrt{1-4z+2z^2-4z^3+z^4}}{2z^3(1-z+z^2)},\\
f(z)&=\frac{1-2z-z^2-\sqrt{1-4z+2z^2-4z^3+z^4}}{2z}.\\
\end{align*}
We let $G(z,u)=g_0(z)+zu\frac{g(z)}{1-uf(z)}$, the bivariate generating function of the almost Riordan array of first order. We seek to find $G(z,\frac{z}{1-z})$. We get
$$G(z,\frac{z}{1-z})=\frac{(1-u+zu)(1-2z+2zu-(1+2u)z^2-\sqrt{1-4z+2z^2-4z^3+z^4}}{2(1-z+z^2)(u(u-1)+z(1-u^2)+z^2u(u+1))},$$ which coincides with the generating function  $\mathit{Total}(z,u)$ of the matrix $\mathcal{T}$.
\hfill $\Box$
\medskip

\begin{prop} The matrix \[\left(
\begin{array}{cccccc}
  1 & 1 & 1 & 1 & 1 &\cdots  \\
  3 & 3 & 3 & 3 & 3 &\cdots  \\
  8 & 10 & 12 & 14 & 16  &\cdots  \\
  23 & 33 & 43 & 53 & 63 &\cdots  \\
  69 & 107 & 149 & 195 & 245 &\cdots \\
  \vdots & \vdots & \vdots & \vdots & \vdots &\ddots \\
\end{array}
\right)\] is the rectification of the Riordan array $(g(z),f(z))$ with 
\begin{align*}f(z)&=\frac{1-z^2-\sqrt{1-4z+2z^2-4z^3+z^4}}{2},\quad\mbox{and}\\
g(z)&=\frac{1-3z+z^2-z^3-(1-z)\sqrt{1-4z+2z^2-4z^3+z^4}}{2z^3(1-z+z^2)}.
\end{align*}
\end{prop}
\noindent {\it Proof.} It suffices to check that the g.f. of $\mathcal{T}$, i.e. $\mathit{Total}(z,u)$, equals to $g_0(z)+zu\frac{g(z)}{1-u\frac{f(z)}{z}}$.
\hfill $\Box$

We can express $f(z)$ and $g(z)$, respectively, in the following form
$$g(z)=\frac{1}{1-3z+z^2-z^3}C\left(\frac{z^3(1-z+z^2)}{(1-3z+z^2-z^3)^2}\right),$$
$$f(z)=\frac{z(1-z+z^2)}{1-z^2}C\left(\frac{z(1-z+z^2)}{(1-z^2)^2}\right).$$
Then $g(z)$ expands to give the first column $1,3,8,23,\ldots$, whose $n$-th term $v_n$ can be expressed
$$v_n=\scalemath{0.75}{\sum\limits_{k=0}^{n}\sum\limits_{j=0}^{k}\binom{k}{j}(-1)^j\sum\limits_{i=0}^{j}\binom{j}{i}(-1)^i\sum\limits_{\ell=0}^{n-3k-j-i}\binom{2k+\ell}{\ell}\sum\limits_{m=0}^{\ell}\binom{\ell}{m}3^{\ell-m}(-1)^m\binom{m}{n-3k-j-i-\ell-m}(-1)^{n-3k-j-i-\ell-m}c_k}.$$
Using $v_n$, we can deduce the following.
\begin{prop}The general term $t(n,k)$ equals
$$\scalemath{0.75}{\sum\limits_{i=0}^{n+k}v_{n+k-i}\sum\limits_{j=0}^{i}\sum\limits_{m=0}^{j}M_{m,k}\sum\limits_{p=0}^{m}\binom{m}{p}(-1)^p\sum\limits_{q=0}^{p}\binom{p}{q}(-1)^q\binom{2m-1+\frac{j-m-p-q}{2}}{\frac{j-m-p-q}{2}}\frac{1+(-1)^{j-m-p-q}}{2}\binom{k}{\frac{i-j}{2}}(-1)^{\frac{i-j}{2}}\frac{1+(-1)^{i-j}}{2}},  $$
where 
$$M_{n,k}=\left\{\begin{array}{ll}
[k=0]& \mbox{ if } n=0\\
\frac{n}{k}\binom{2n-k-1}{n-k}&\mbox{ otherwise}\end{array}\right.$$
is the general term of Riordan array $(1,zC(z))$ (see \oeis{A106566}).
\end{prop}

\subsection{Comment on Section 3.1}
\begin{prop}
The matrix $\mathcal{T}=[t(n,k)]_{n\ge0, k \ge0}$  is the Riordan array
\begin{align*}\mathcal{T}&=\left(1+zM(z),z(1+zM(z))\right)\\
&=\left(C\left(\frac{z}{1+z}\right),zC\left(\frac{z}{1+z}\right)\right),
\end{align*}
where $C(z)=\frac{1-\sqrt{1-4z}}{2z}$ is the generating function of the Catalan numbers $c_n=\frac{1}{n+1}\binom{2n}{n}$, and $M(z)=\frac{1-z-\sqrt{1-2z-3z^2}}{2z^2}$ is the generating function of the Motzkin numbers \oeis{A001006}.
\end{prop}
\noindent {\it Proof.} It suffices to check that  $\mathit{Total}(z,u)=\frac{g(z)}{1-uf(z)}$.
\hfill $\Box$
\medskip

This triangle $\mathcal{T}$ corresponds to  \oeis{A091836} in \cite{oeis} where the coefficient of row $n-1$ and column $k$ is the number of Motzkin paths of length $n$ having $k$ points on the horizontal axis (besides the first and last point). 

As a consequence, we deduce that 
$$t(n,k)=\left\{\begin{array}{ll}1, & \mbox{if } n=k,\\
\frac{k+1}{n+1}\sum\limits_{j=1}^{n-k} j(-1)^{n-k-j}\binom{n+j}{j}\sum\limits_{i=0}^{n-k}\frac{1}{n-k}\binom{i}{n-k-i+j}\binom{n-k}{i},& \mbox{otherwise.}\end{array}\right.
$$
Alternatively, we have 
$$t(n,k)=\sum_{j=0}^{n}\binom{n-1}{n-j}(-1)^{n-j}(k+1)\sum_{i=0}^{j-k}\frac{1}{k+i+1}\binom{k+2i}{i}\binom{j-i-1}{j-k-i}.$$
A third expression for $t(n,k)$ is given by the following proposition.
\begin{prop} The general term $t(n,k)$ of the Riordan array $(1 + zM(z), z(1 + zM(z)))$ is
given by
$$t(n,k)=\left\{\begin{array}{ll}
1& \mbox{ if }n=k\\
\frac{k+1}{n-k}\sum_{j=0}^k\binom{k}{j}\sum_{i=0}^{n-k}\binom{n-k}{i}\binom{i}{n-k-i-j-1}&\mbox{ otherwise}
\end{array}\right..
$$
\end{prop}
\noindent {\it Proof.} We prove this using Lagrange inversion, using the fact that $$(zM(z))^{(-1)}=\frac{z}{1+z+z^2}.$$
Thus we have
\begin{align*} t(n,k)&=[z^n](1+zM(z))(z(1+zM(z)))^k\\
&=[z^{n-k}](1+zM(z))^{k+1}\\
&=[z^{n-k}]G(zM(z)), \quad \mbox{ with }\quad  G(z)=(1+z)^{k+1}\\
&=\frac{1}{n-k}[z^{n-k-1}]G'(z)\left(\frac{z}{(zM(z))^{(-1)}}\right)^{n-k} \mbox{ (Lagrange inversion) }\\
&=\frac{1}{n-k}[z^{n-k-1}](k+1)(1+z)^k(1+z+z^2)^{n-k}\\
&=\frac{k+1}{n-k} [z^{n-k}]\sum\limits_{j=0}^k\binom{k}{j}z^j\sum\limits_{i=0}^{n-k}\binom{n-k}{i}z^i(1+z)^i\\
&=\frac{k+1}{n-k} [z^{n-k}]\sum\limits_{j=0}^k\binom{k}{j}z^j\sum\limits_{i=0}^{n-k}\binom{n-k}{i}z^i\sum\limits_{\ell=0}^{i}\binom{i}{\ell}z^\ell\\
&=\frac{k+1}{n-k}\sum_{j=0}^k\binom{k}{j}\sum_{i=0}^{n-k}\binom{n-k}{i}\binom{i}{n-k-i-j-1}.
\end{align*}
\hfill $\Box$

\subsection{Comment on Section 3.2}

\begin{prop} The matrix $\mathcal{T}=[t(n,k)]_{n\geq 0, k \geq 0}$ can be written  
\[\mathcal{T}=
\left(
\begin{array}{ccccccc}
 1 & 0 & 0 & 0 & 0 & 0 & \cdots  \\
 1 & 1 & 1 & 1 & 1 & 1 & \cdots \\
 1 & 1 & 1 & 1 & 1 & 1 & \cdots \\
 2 & 3 & 4 & 5 & 6 & 7 & \cdots \\
 4 & 6 & 8 & 10 & 12 & 14 & \cdots\\
 9 & 15 & 22 & 30 & 39 & 49 & \cdots\\
     \vdots & \vdots & \vdots & \vdots & \vdots& \vdots  &\ddots \\
\end{array}
\right)=A\cdot B\]
where \[ A=\left(
\begin{array}{ccccccc}
 1 & 0 & 0 & 0 & 0 & 0 & \cdots \\
 1 & 1 & 0 & 0 & 0 & 0 & \cdots \\
 1 & 1 & 0 & 0 & 0 & 0 & \cdots \\
 2 & 3 & 1 & 0 & 0 & 0 & \cdots \\
 4 & 6 & 2 & 0 & 0 & 0 & \cdots \\
 9 & 15 & 7 & 1 & 0 & 0 & \cdots\\
  \vdots & \vdots & \vdots & \vdots & \vdots& \vdots  &\ddots \\
\end{array}
\right)\quad \mbox{and}\quad B=\left(
\begin{array}{ccccccc}
 1 & 0 & 0 & 0 & 0 & 0 & \cdots \\
 0 & 1 & 1 & 1 & 1 & 1 & \cdots \\
 0 & 0 & 1 & 2 & 3 & 4 & \cdots \\
 0 & 0 & 0 & 1 & 3 & 6 & \cdots \\
 0 & 0 & 0 & 0 & 1 & 4 & \cdots \\
 0 & 0 & 0 & 0 & 0 & 1 & \cdots\\
  \vdots & \vdots & \vdots & \vdots & \vdots& \vdots  &\ddots \\
\end{array}
\right)
\]
are defined as follows:

$\bullet$ The matrix $B=[b_{n,k}]_{n,k\geq 0}$ is defined by $b_{0,0}=1$, and $b_{n,0}=b_{0,n}=0$ if $n\geq 1$, and $b_{n,k}=\binom{k-1}{n-1}$ otherwise, which is a kind of Pascal matrix.

$\bullet$ The matrix $A=[a_{n,k}]_{n,k\geq 0}$ is the almost `stretched' Riordan array with initial column whose generating function is 
$$g_0(z)=1+zM(z)=\frac{1+z-\sqrt{1-2z-3z^2}}{2z},$$
which is followed by the shifted stretched Riordan array $(g(z),z^2g(z))$ where $$g(z)=\frac{1-z-2z^2-\sqrt{1-2z-3z^2}}{2z^3(1+z)}.$$
\end{prop}
\noindent {\it Proof.} The almost Riordan array $A$ has generating function 
$$g_0(z)+zu\frac{g(z)}{1-z^2ug(z)}.$$
Using the fundamental theorem of Riordan arrays, the product has generating function 
$$g_0(z)+x\frac{u}{1-u}\frac{g(z)}{1-z^2\frac{u}{1-u}g(z)}.$$
By simplifying this expression, we obtain $\mathit{Total}(z,u)$, which completes the proof.
\hfill $\Box$

\begin{prop} The matrix \[\left(
\begin{array}{cccccc}
  1 & 1 & 1 & 1 & 1 &\cdots  \\
  1 & 1 & 1 & 1 & 1 &\cdots  \\
  3 & 4 & 5 & 6 & 7  &\cdots  \\
  6 & 8 & 10 & 12 & 14 &\cdots  \\
  15 & 22 & 30 & 39 & 49 &\cdots \\
  36 & 54 & 75 & 99 & 126 &\cdots \\
  \vdots & \vdots & \vdots & \vdots & \vdots &\ddots \\
\end{array}
\right)\] is the rectification of the Riordan array $(M(z),zR(z))$ where $M(z)=\frac{1-z-\sqrt{1-2z-3z^2}}{2z^2}$ is the g.f. of the Motzkin numbers, and $R(z)=\frac{1+z-\sqrt{1-2z-3z^2}}{2z(1+z)}$ is the g.f. of the Riordan numbers (\oeis{A005043}).
\end{prop}
\noindent {\it Proof.} It suffices to check that the generating function of $\mathcal{T}$, i.e. $\mathit{Total}(z,u)$, equals to   $g_0(z)+zu\frac{M(z)}{1-u\frac{zR(z)}{z}}$.
\hfill $\Box$

We let $m_n$ denote the $n$-th Motzkin number $m_n=\sum_{k=0}^{\lfloor\frac{n}{2}\rfloor}\binom{n}{2k}c_k$ where $c_k$ is the $k$-th Catalan defined above.
\begin{cor} We have 
$$t(n,k)=\left\{
\begin{array}{ll} [k=0]& \mbox{ if } n=0\\
r(n-1,k)&\mbox{ otherwise }
\end{array}
\right., $$
where 
$$r(n,k)=\scalemath{0.85}{\sum\limits_{i=0}^n\frac{m_i\cdot (k+[n=i])}{n+k-i+[n=i]}\sum\limits_{j=0}^{n-i}(-1)^{n-i-j}\binom{n+k-i+j-1}{j}\sum\limits_{\ell=0}^{n+k-i}\binom{n+k-i}{\ell}\binom{\ell}{n-i-j-\ell}}.$$
\end{cor}
\noindent {\it Proof.}
We have $(M(z),zR(z))^{-1}=\left(\frac{(1-z)^2}{1-z+z^2},\frac{z(1-z)}{1-z+z^2}\right).$ If we denote by $(v(z),u(z))$ these inverse Riordan array, then we obtain 
$$(M(z),zR(z))=\left( \frac{1}{v(\bar{u}(z))},\bar{u}(z)\right).$$
Using the definition of a Riordan array, and Lagrange inversion, we find that the Riordan array $(M(z),zR(z))$ has general term $r(n,k)$ given by
$$r(n,k)=\scalemath{0.85}{\sum\limits_{i=0}^n\frac{m_i\cdot (k+[n=k+i])}{n-i+[n=k+i]}\sum\limits_{j=0}^{n-k-i}(-1)^{n-k-i-j}\binom{n-i+j-1}{j}\sum\limits_{\ell=0}^{n-i}\binom{n-i}{\ell}\binom{\ell}{n-k-i-j-\ell}}.$$
To rectify this array, we change $n$ to $n+k$, and the result follows.
\hfill $\Box$

\begin{rem} The Riordan array $(1+zM(z), z(1+zM(z)))$ is a pseudo-involution in the Riordan group (see \cite{Bur}), that is, the matrix $[(-1)^k t_{n,k}]_{n,k\geq 0}$ is idempotent. Thus, this work yields a significant lattice path interpretation of this array.
\end{rem}

\end{document}